\documentclass[12 pt]{amsart}
\usepackage{graphicx}

\usepackage[pagebackref]{hyperref}
\hypersetup{
colorlinks = true,
linkcolor = blue,
citecolor = blue}
\usepackage[utf8]{inputenc}
\usepackage[T1]{fontenc}

\usepackage[]{amsfonts}
\usepackage{amssymb}
\usepackage{amsmath}
\def\couleur(#1 #2 #3)
	{
\special{ps::[end]}}

\def\bx#1{\setbox1=\hbox{\kern3pt{#1}\kern3pt}			
 \dimen1=\ht1 \advance\dimen1 by 3pt \dimen2=\dp1 \advance\dimen2 by 3pt
 \setbox1=\hbox{\vrule height\dimen1 depth\dimen2\box1\vrule}%
 \setbox1=\vbox{\hrule\box1\hrule}%
 \advance\dimen1 by .4pt \ht1=\dimen1
 \advance\dimen2 by .4pt \dp1=\dimen2 \box1\relax}

\def\wbb#1{\kern#1em}
\def\vci{\vrule  width.02em height1.47ex depth-.0ex}		
\def\11{{\rm\wbb{.2}\vci\wbb{-.37}1}}

\def\Supp{\mathop{\rm Supp}\nolimits}
\def\underset#1#2{\mathrel{\mathop{\kern0pt #2}\limits_{#1}}}

\def\overset#1#2{\mathrel{\mathop{\kern0pt #2}\limits^{#1}}}

\newtheorem{thm}{Theorem}[section]
\newtheorem{lem}[thm]{Lemma}
\newtheorem{prop}[thm]{Proposition}
\newtheorem{cor}[thm]{Corollary}
\newtheorem{defin}[thm]{Definition}
\newtheorem{rem}[thm]{Remark}
\newtheorem{exa}[thm]{Example}

\begin{document}
\renewcommand{\abstractname}{Abstract}

\begin{abstract}
We study the $\bar \partial $-equation first in Stein manifold
 then in complete K\"ahler manifolds. The aim is to get $L^{r}$
 and Sobolev estimates on solutions with compact support.\par 
\quad In the Stein case we get that for any $(p,q)$-form $\omega $
 in $L^{r}$ with compact support and $\bar \partial $-closed
 there is a $(p,q-1)$-form $u$ in $W^{1,r}$ with compact support
 and such that $\bar \partial u=\omega .$\par 
\quad In the case of K\"ahler manifold, we prove and use estimates
 on solutions on Poisson equation with compact support and the
 link with $\bar \partial $ equation is done by a classical theorem
 stating that the Hodge laplacian is twice the $\bar \partial
 $ (or Kohn) Laplacian in a K\"ahler manifold.\par 
\quad This uses and improves, in special cases, our result on Andreotti-Grauert
 type theorem.\par 

\end{abstract}

\title{Solutions of the $\bar \partial $-equation with compact support on Stein and K\"ahler manifold.}

\author{Eric Amar}

\date{}

\keywords{$\bar \partial $-equation, Poisson equation, $L^{r}$ estimates, Stein, Rieman and K\"ahler manifolds.}
\maketitle
 \ \par 
\ \par 
\ \par 
\ \par 
\ \par 
\ \par 

\tableofcontents
\ \par 

\section{Introduction.}
\quad The study of $L^{r}$ solutions for the $\bar \partial $ equation
 is a long standing subject. By use of PDE methods, L. H\"ormander~\cite{Hormander73}
 get his famous $L^{2}$ estimates, we shall still use here.\ \par 
\quad The next results were obtained by the use of solving kernels:
 Grauert-Lieb~\cite{GrauLieb70}, Henkin~\cite{Henkin70}, Ovrelid~\cite{Ovrelid71},
 Skoda~\cite{zeroSkoda}, Krantz~\cite{KrantzDbar76}, in the case
 of  strictly pseudo-convex domains with $\displaystyle {\mathcal{C}}^{\infty
 }$ smooth boundary in $\displaystyle {\mathbb{C}}^{n},$  with
 the exception of  Kerzman~\cite{Kerzman71} who worked in the
 case of $\displaystyle (0,1)$  forms in  strictly pseudo-convex
 domains with $\displaystyle {\mathcal{C}}^{4}$ smooth boundary
 in Stein manifolds.\ \par 
\quad The case of smooth solutions with compact support goes to the
 work of Andreotti and\ \par 
Grauert~\cite{AndreGrauert62}. Our aim here is to study the same
 problem but with $L^{r}$ estimates, as we did in~\cite{AnGrauLr19}.\ \par 
\quad In a first part, we recall some results on solutions of the $\bar
 \partial $ equation in complex manifolds. Then, modifying a
 result by C. Laurent-Thi\'ebaut~\cite{Claurent15}, we prove that:\ \par 

\begin{cor}
Let $X$ be a complex manifold. Let $S$ be a $(p,q-1)$-current
 in $L_{p,q-1}^{r}(X)$ with compact support $W$ in $X.$ Suppose
 that $W\subset \Omega _{1}\Subset \Omega _{2},$ where $\Omega
 _{j},\ j=1,2$ are relatively compact  pseudo-convex open sets
 with smooth ${\mathcal{C}}^{\infty }$ boundary in $X$ and such
 that  there is a strictly pluri-subharmonic function  $\rho
 _{1}$ in ${\mathcal{C}}^{3}(\bar \Omega _{1}).$\par 
Moreover suppose that $\omega :=\bar \partial S$ is also in $L_{p,q}^{r}(X).$
 Let $U$ be any neighborhood of $W$ contained in $\Omega _{1}.$
 Then there is a $(p,q-1)$-current $u$ with compact support in
 $U$ such that $\bar \partial u=\omega $ and $u\in W_{p,q-1}^{1,r}(U).$
\end{cor}
\quad Then we get the following result which seems to end a question
 Guiseppe Tomassini ask me almost ten years ago (see ~\cite{AmMong12}
 and ~\cite{AnGrauLr12}).\ \par 
\quad Let $M$ be a complex manifold and $\Lambda _{p,q}(\bar M)$ the
 set of $(p,q)$-forms ${\mathcal{C}}^{\infty }$ in $\bar M.$\ \par 
\quad Recall that the Kohn laplacian $\Delta _{\bar \partial }$ is defined as:\ \par 
\quad \quad \quad $\forall \alpha \in L_{p,q}^{2}(M),\ \Delta _{\bar \partial }\alpha
 :=\bar \partial \bar \partial ^{*}\alpha +\bar \partial ^{*}\bar
 \partial \alpha .$\ \par 
First let us define, as in p. 278 in~\cite{Kohn73}, the harmonic fields:\ \par 
\quad \quad \quad ${\mathcal{H}}_{p,q}:=\lbrace h\in \Lambda _{p,q}(\bar M)::\bar
 \partial h=\bar \partial ^{*}h=0\rbrace .$\ \par 
Then we have:\ \par 

\begin{thm}
~\label{k56} Let $X$ be a Stein manifold and  $\omega $ be a
 $(p,q)$ form in $L^{r}(X),\ r>1$ with compact support in $X.$
 Suppose that $\omega $ is such that:\par 
\quad $\bullet $ if $1\leq q<n,\ \bar \partial \omega =0$; \par 
\quad $\bullet $ if $q=n,\ \forall V\subset X,\ \Supp \omega \subset
 V,\ \omega \perp {\mathcal{H}}_{n-p,0}(V)$\!\!\!\! .\par 
Then there is a $(p,q-1)$ form $u$ in $W^{1,r}(X)$ with compact
 support in $X$ such that $\bar \partial u=\omega $ as distributions
 and 	$\displaystyle {\left\Vert{u}\right\Vert}_{W^{1,r}(\Omega
 )}\leq C{\left\Vert{\omega }\right\Vert}_{L^{r}(\Omega )}.$
\end{thm}
\ \par 
\quad In a second part we study this problem in a K\"ahler manifold.
 The method is completely different:  we first study  $L^{r}$
 solutions with compact support for the Poisson equation in a
 riemannian manifold and we use the link done by the following
 classical theorem relying the Hodge laplacian and the $\bar
 \partial $ (or Kohn) laplacian. See for instance C. Voisin's
 book~\cite{Voisin02}.\ \par 

\begin{thm}
~\label{k60}Let $(X,\kappa )$ be a k\"ahlerian manifold. Let
 $\Delta ,\ \Delta _{\partial },\ \Delta _{\bar \partial }$ the
 laplacians associated to $d,\ \partial ,\ \bar \partial $ respectively.
 Then we have the relations:\par 
\quad \quad \quad $\Delta =2\Delta _{\partial }=2\Delta _{\bar \partial }.$
\end{thm}
\quad We get, with this time ${\mathcal{H}}_{q}(\Omega )=\lbrace h\in
 \Lambda _{q}(\bar \Omega )::\Delta h=0\rbrace $:\ \par 

\begin{thm}
~\label{kF39}Let $(X,\omega )$ be a complete k\"ahlerian manifold.
 Let $\Omega $ be a relatively compact  domain in $X.$ Let $\displaystyle
 \omega \in L_{p,q}^{r}(\Omega ),\ \bar \partial \omega =0$ in
 $\Omega $ and $\omega $ compactly supported in $\Omega .$ Suppose
 moreover that $\displaystyle \omega \perp {\mathcal{H}}_{2n-p-q}(\Omega
 ).$\par 
\quad Then there is a $\displaystyle u\in W_{p,q-1}^{1,r}(\Omega )$
 with compact support in $\Omega $ and such that $\bar \partial u=\omega .$
\end{thm}
\quad This result seems weaker than the previous one because we need
 that $\displaystyle \omega \perp {\mathcal{H}}_{2n-p-q}(\Omega
 ),$ but, unless $X$ is weakly pseudo-convex, a compact set is
 not contained in a pseudo-convex one in general. Hence the method
 used for the proof of Theorem~\ref{k56} cannot apply here.\ \par 
\ \par 
\quad This work is presented the following way.\ \par 
\ \par 
For the first part:\ \par 
\quad $\bullet $ In Section~\ref{k57} we recall results on estimates
 for the $\bar \partial $ equation.\ \par 
\ \par 
\quad $\bullet $ In Subsection~\ref{k58} we recall the notion of $r$-regularity
 and its consequence in term of solution of the $\bar \partial
 $ equation with compact support.\ \par 
\ \par 
\quad $\bullet $ In Subsection~\ref{k59} we show that, under some circumstances,
 the regularity of solutions of the $\bar \partial $ equation
 may increases.\ \par 
This part is directly coming from a work of C. Laurent-Thi\'ebaut~\cite{Claurent15}.\
 \par 
\ \par 
For the second part:\ \par 
\quad $\bullet $ In Section~\ref{kF40} we start with the Hodge laplacian
 on a riemannian manifold and we recall results we get in~\cite{ellipticEq18}
 concerning the Poisson equation.\ \par 
\ \par 
\quad $\bullet $ In Section~\ref{kF42} we study the solutions of the
 Poisson equation with compact support and we prove, using weighted
 estimates:\ \par 

\begin{thm}
Let $X$ be a complete oriented riemannian manifold. Let $\Omega
 $ be a relatively compact domain in $X.$ Let $\omega \in L^{r}_{p}(\Omega
 )$ with compact support in $\Omega $ and such that $\omega $
 is orthogonal to the harmonic $p$-forms $\displaystyle {\mathcal{H}}_{p}(\Omega
 ).$ Then there is a $p$-form $u\in W^{2,r}_{p}(\Omega )$ with
 compact support in $\Omega $ such that $\Delta u=\omega $ as
 distributions and 	$\displaystyle {\left\Vert{u}\right\Vert}_{W^{2,r}_{p}(\Omega
 )}\leq C{\left\Vert{\omega }\right\Vert}_{L^{r}_{p}(\Omega )}.$
\end{thm}
\quad $\bullet $ In Section~\ref{kF41}, using equality of the laplacians,
 we prove Theorem~\ref{kF39}.\ \par 
\ \par 
\quad $\bullet $ Finally in the Appendix we prove certainly known results
 on the duality $L^{r}-L^{r'}$ for $(p,q)$-forms in a complex
 manifold we needed.\ \par 

\section{On estimates for the $\displaystyle \bar \partial $
 equation in complex manifolds.~\label{k57}}
\quad Here we shall be interested in strictly $c$-convex (s.c.c. for
 short) domains $D$ in a complex manifold. Such a domain is defined
 by a function  $\rho $ of class $\displaystyle {\mathcal{C}}^{3}$
 in a neighbourhood $U$ of $\bar D$  and such that  $\displaystyle
 i\partial \bar \partial \rho $  has at least $\displaystyle
 n-c+1$  strictly positive eigenvalues in $U.$\ \par 
\ \par 
\quad We have the following Theorem 1.1 from~\cite{dBarIntersection17}:\ \par 

\begin{thm}
~\label{k49}Let $\Omega $ be a Stein manifold of dimension $n$
 and a s.c.c. domain $D$ such that $D$ is relatively compact
 with smooth $\displaystyle {\mathcal{C}}^{3}$ boundary in $\Omega
 .$ Let  $\omega $ be a $\displaystyle (p,q)$  form in $\displaystyle
 L^{r}_{p,q}(D),\ \bar \partial \omega =0$  with $\displaystyle
 1<r<2n+2,\ c\leq q\leq n.$ Then there is a $\displaystyle (p,q-1)$
  form $u$ in $\displaystyle L^{s}(D),$ with $\displaystyle \
 \frac{1}{s}=\frac{1}{r}-\frac{1}{2n+2},$ such that $\displaystyle
 \bar \partial u=\omega .$\par 
\quad If  $\omega $ is in $\displaystyle L^{r}_{p,q}(D),\ \bar \partial
 \omega =0$  with  $\displaystyle r\geq 2n+2,\ c\leq q\leq n,$
 then there is a $\displaystyle (p,q-1)$  form $u$ in $\displaystyle
 \Lambda _{(p,q-1)}^{\epsilon }(\bar D)$ such that $\displaystyle
 \bar \partial u=\omega $ with $\displaystyle \epsilon =\frac{1}{2}-\frac{n+1}{r}.$
\end{thm}
The
 spaces $\displaystyle \Lambda _{(p,q-1)}^{\epsilon }(\bar D)$
 are the (isotropic) Lipschitz spaces of order $\epsilon $ and
 we set $\displaystyle \ \Lambda _{(p,q-1)}^{0}(\bar D):=L_{(p,q-1)}^{\infty
 }(D).$\ \par 
\quad This theorem has the obvious corollary:\ \par 

\begin{cor}
~\label{k51}Let $\Omega $ be a complex manifold of dimension
 $n$ and a domain $D$ relatively compact with smooth $\displaystyle
 {\mathcal{C}}^{3}$ boundary in $\Omega .$ Suppose moreover that
 $D:=\lbrace \rho <0\rbrace ,$ where $\rho $ is a strictly pluri-subharmonic
 function in ${\mathcal{C}}^{\infty }(\bar D)$ with $\left\vert{\partial
 \rho }\right\vert >0$ on $\partial D.$ Let  $\omega $ be a $\displaystyle
 (p,q)$  form in $\displaystyle L^{r}_{p,q}(D),\ \bar \partial
 \omega =0$  with $\displaystyle 1<r<2n+2,\ 1\leq q\leq n.$ Then
 there is a $\displaystyle (p,q-1)$  form $u$ in $\displaystyle
 L^{s}(D),$ with $\displaystyle \ \frac{1}{s}=\frac{1}{r}-\frac{1}{2n+2},$
 such that $\displaystyle \bar \partial u=\omega ,$ with  ${\left\Vert{u}\right\Vert}_{s}\leq
 C{\left\Vert{\omega }\right\Vert}_{r}.$\par 
\quad If  $\omega $ is in $\displaystyle L^{r}_{p,q}(D),\ \bar \partial
 \omega =0$  with  $\displaystyle r\geq 2n+2,\ 1\leq q\leq n,$
 then there is a $\displaystyle (p,q-1)$  form $u$ in $\displaystyle
 \Lambda _{(p,q-1)}^{\epsilon }(\bar D)$ such that $\displaystyle
 \bar \partial u=\omega $ with $\displaystyle \epsilon =\frac{1}{2}-\frac{n+1}{r}$
 and  ${\left\Vert{u}\right\Vert}_{\Lambda ^{\epsilon }}\leq
 C{\left\Vert{\omega }\right\Vert}_{r}.$
\end{cor}
\quad Proof.\ \par 
Take a convex increasing function $\chi $ on ${\mathbb{R}}^{-}$
 such that $\displaystyle \chi (t)\rightarrow \infty $ when $t\rightarrow
 0.$ The function $\varphi (z):=\chi \circ \rho (z)$ is still
 strictly pluri-subharmonic on $D$ and exhausting. So $D$ is
 a Stein manifold by Theorem 5.2.10 in~\cite{Hormander73}. A
 strictly pseudo-convex domain is a s.c.c. domain with $c=1,$
 so we can apply Theorem~\ref{k49}. $\hfill\blacksquare $\ \par 
\ \par 

\begin{cor}
~\label{k52}Let $\Omega $ be a complex manifold of dimension
 $n$ and a domain $D$ relatively compact with smooth $\displaystyle
 {\mathcal{C}}^{3}$ boundary in $\Omega .$ Suppose moreover that
 $D:=\lbrace \rho <0\rbrace ,$ where $\rho $ is a strictly pluri-subharmonic
 function in ${\mathcal{C}}^{3}(\bar D)$ with $\left\vert{\partial
 \rho }\right\vert >0$ on $\partial D.$ Let  $\omega $ be a $\displaystyle
 (p,q)$  form in $\displaystyle L^{r}_{p,q}(D),\ \bar \partial
 \omega =0$  with $\displaystyle 1<r<\infty ,\ 1\leq q\leq n.$
 Then there is a $\displaystyle (p,q-1)$  form $u$ in $\displaystyle
 L^{r}(D),$ with  ${\left\Vert{u}\right\Vert}_{r}\leq C{\left\Vert{\omega
 }\right\Vert}_{r},$ such that $\displaystyle \bar \partial u=\omega .$
\end{cor}
\quad Proof.\ \par 
Because $D$ is relatively compact, if $\displaystyle u\in L^{s}(D)$
 for $s\geq r$ then $\displaystyle u\in L^{r}(D).$ Then the Corollary~\ref{k51}
 gives the result. $\hfill\blacksquare $\ \par 

\subsection{Weak $r$-regularity.~\label{k58}}
\ \par 
\quad We shall need the definition, see~\cite{AnGrauLr19}:\ \par 

\begin{defin}
~\label{AG16}Let $X$ be a complex manifold equipped with a Borel
 $\sigma $-finite measure $dm$ and $\Omega $ a domain in $X$;
 let $\displaystyle r\in \lbrack 1,\ \infty \rbrack ,$ we shall
 say that $\Omega $ is $r$ {\bf regular} if for any $p,q\in \lbrace
 0,...,n\rbrace ,\ q\geq 1,$ there is a constant $C=C_{p,q}(\Omega
 )$ such that for any $(p,q)$ form $\omega ,\ \bar \partial $
 closed in $\Omega $ and in $L^{r}(\Omega ,dm)$ there is a $(p,q-1)$
 form $u\in L^{r}(\Omega ,dm)$ such that $\bar \partial u=\omega
 $ and $\ {\left\Vert{u}\right\Vert}_{L^{r}(\Omega )}\leq C{\left\Vert{\omega
 }\right\Vert}_{L^{r}(\Omega )}.$\par 
\quad \quad 	We shall say that $\Omega $ is {\bf weakly }$r$ {\bf regular}
 if for any compact set $K\Subset \Omega $ there are $3$ open
 sets $\Omega _{1},\Omega _{2},\Omega _{3}$ such that $K\Subset
 \Omega _{3}\subset \Omega _{2}\subset \Omega _{1}\subset \Omega
 _{0}:=\Omega $ and $3$ constants $C_{1},C_{2},C_{3}$ such that:\par 
\quad \quad \quad $\displaystyle \forall j=0,1,2,\ \forall p,q\in \lbrace 0,...,n\rbrace
 ,\ q\geq 1,\ \forall \omega \in L_{p,q}^{r}(\Omega _{j},dm),\
 \bar \partial \omega =0,$\par 
\quad \quad \quad \quad \quad $\displaystyle \exists u\in L_{p,q-1}^{r}(\Omega _{j+1},dm),\
 \bar \partial u=\omega $\par 
and $\ {\left\Vert{u}\right\Vert}_{L^{r}(\Omega _{j+1})}\leq
 C_{j+1}{\left\Vert{\omega }\right\Vert}_{L^{r}(\Omega _{j})}.$\par 
\quad \quad 	I.e. we have a $3$ steps chain of resolution.
\end{defin}
\ \par 
\quad in~\cite{AnGrauLr19} we prove the Theorem 3.5, p. 6, where $L^{r,c}(\Omega
 )$ means that the form is in $L^{r}(\Omega )$ with compact support
 in $\Omega $:\ \par 

\begin{thm}
~\label{k46}Let $\Omega $ be a weakly $r'$ regular domain in
 a complex manifold and  $\omega $ be a $(p,q)$ form in $L^{r,c}(\Omega
 ),\ r>1.$ Suppose that $\omega $ is such that:\par 
\quad $\bullet $ if $1\leq q<n,\ \bar \partial \omega =0$; \par 
\quad $\bullet $ if $q=n,\ \forall V\subset \Omega ,\ \Supp \omega
 \subset V,\ \omega \perp {\mathcal{H}}_{n-p}(V)$\!\!\!\! .\par 
Then there is a $C>0$ and a $(p,q-1)$ form $u$ in $L^{r,c}(\Omega
 )$ such that $\bar \partial u=\omega $ as distributions and
 	$\displaystyle {\left\Vert{u}\right\Vert}_{L^{r}(\Omega )}\leq
 C{\left\Vert{\omega }\right\Vert}_{L^{r}(\Omega )}.$
\end{thm}
\quad In fact in~\cite{AnGrauLr19} we made the general assumption that
 our complex manifold $X$ is Stein, just to be sure that \emph{any}
 compact set is in a weakly regular domain, because we proved
 in~\cite{AnGrauLr19} that a Stein manifold is weakly $r$-regular.
 The example~\ref{SI0} prove that this is not the case in general.
 But in Theorem~\ref{k46}, the proof works for $X$ being just
 a complex manifold.\ \par 
\quad As a corollary we get:\ \par 

\begin{cor}
~\label{k53} Let $\Omega $ be a complex manifold of dimension
 $n$ and a domain $D$ relatively compact with smooth $\displaystyle
 {\mathcal{C}}^{3}$ boundary in $\Omega .$ Suppose moreover that
 $D:=\lbrace \rho <0\rbrace ,$ where $\rho $ is a strictly pluri-subharmonic
 function in ${\mathcal{C}}^{3}(\bar D)$ with $\left\vert{\partial
 \rho }\right\vert >0$ on $\partial D.$ Then $D$ is $r'$-regular.
 Moreover suppose that  $\omega $ is a $(p,q)$ form in $L^{r,c}(D),\
 r>1$ such that:\par 
\quad $\bullet $ if $1\leq q<n,\ \bar \partial \omega =0$; \par 
\quad $\bullet $ if $q=n,\ \forall V\subset D,\ \Supp \omega \subset
 V,\ \omega \perp {\mathcal{H}}_{n-p}(V)$\!\!\!\! .\par 
Then there is a $C>0$ and a $(p,q-1)$ form $u$ in $L^{r,c}(D)$
 such that $\bar \partial u=\omega $ as distributions and 	$\displaystyle
 {\left\Vert{u}\right\Vert}_{L^{r}(D)}\leq C{\left\Vert{\omega
 }\right\Vert}_{L^{r}(D)}.$
\end{cor}
\quad Proof.\ \par 
The Corollary~\ref{k52} gives that $D$ is $r'$-regular for any
 $1<r'<\infty .$ Hence we can apply Theorem~\ref{k46} to the
 domain $D.$ $\hfill\blacksquare $\ \par 

\subsection{On an improvement of regularity.~\label{k59}}
\ \par 
\quad This section is coming from C. Laurent-Thi\'ebaut~\cite{Claurent15},
 Proposition 1.4 p. 257.\ \par 
To use Theorem 4 in~\cite{BealsGreiStan87} we need a compact
 complex manifold with a smooth ${\mathcal{C}}^{\infty }$ boundary
 having property $Z(q)$ and a $(p,q)$-form in the range of the
 Kohn laplacian, which means that the form must be orthogonal
 to the harmonic fields.\ \par 
\quad In~\cite{Claurent15} the author used twice this Theorem without
 any references to these  two conditions. Because I was unable
 to understand why they are fulfilled, I prove here a weaker
 result fitting well enough with my purpose.\ \par 
\quad The nice idea to work with \emph{exact forms} the regularity
 of which being increased is due to C. Laurent-Thi\'ebaut.\ \par 
\ \par 
\quad Let us define, as in p. 278 in~\cite{Kohn73}, the harmonic fields:\ \par 
\quad \quad \quad ${\mathcal{H}}_{p,q}:=\lbrace h\in \Lambda _{p,q}(\bar M)::\bar
 \partial h=\bar \partial ^{*}h=0\rbrace .$\ \par 
Then we have in~\cite{Kohn73}, that:\ \par 
\quad $\bullet $ ${\mathcal{H}}_{p,q}=\mathrm{k}\mathrm{e}\mathrm{r}\Delta
 _{\bar \partial };$\ \par 
\quad $\bullet $ ${\mathcal{H}}_{p,q}$ is a finite dimensional subspace
 of $\Lambda _{p,q}(\bar M).$\ \par 
We shall denote $H$ the orthogonal projection of $L_{p,q}^{2}(M)$
 onto ${\mathcal{H}}_{p,q}.$\ \par 
And we have a Hodge decomposition, eq. 2.26, p. 278 in~\cite{Kohn73}:\ \par 
\quad \quad \quad $\forall \alpha \in L_{p,q}^{2}(M),\ \alpha =\bar \partial \bar
 \partial ^{*}\alpha +\bar \partial ^{*}\bar \partial \alpha +H\alpha .$\ \par 
Moreover if $\bar \partial \alpha =0$ and $\alpha \perp {\mathcal{H}}_{n-p,n-q}$
 then, eq. 2.27, p. 278 in~\cite{Kohn73}:\ \par 
\quad \quad \quad $\alpha =\bar \partial \bar \partial ^{*}N\alpha $ and $\varphi
 =\bar \partial ^{*}N\alpha $ is the unique solution of the equation
 $\alpha =\bar \partial \varphi $ orthogonal to ${\mathcal{H}}_{n-p,n+1-q}.$\
 \par 
\quad Recall also that a pseudo-convex set with smooth ${\mathcal{C}}^{\infty
 }$ boundary has the $Z(q)$ property for any $q\geq 1.$ This
 is why we shall use mainly this notion.\ \par 

\begin{thm}
~\label{k45}Let $X$ be a complex manifold and $q\geq 1.$ Let
 $S$ be a $(p,q-1)$-current in $L_{p,q-1}^{r}(X).$  with compact
 support $W$ in $X.$ Suppose that $W\subset \Omega _{1}\Subset
 \Omega _{2},$ where $\Omega _{j},\ j=1,2$ are relatively compact
  pseudo-convex open sets with smooth ${\mathcal{C}}^{\infty
 }$ boundary in $X$ and such that $S\perp {\mathcal{H}}_{n-p,n+1-q}(\Omega
 _{1}).$ \par 
Moreover suppose that $\omega :=\bar \partial S$ is also in $L_{p,q}^{r}(X).$
 Let $U$ be any neighborhood of $W$ contained in $\Omega _{1}.$
 Then there is a $(p,q-1)$-current $u$ with compact support in
 $U$ such that $\bar \partial u=\omega $ and $u\in W_{p,q-1}^{1,r}(U).$
\end{thm}
\quad Proof.\ \par 
As done in~\cite{Claurent15} we shall use Theorem 4 in~\cite{BealsGreiStan87}.\
 \par 
\quad We have to see that $\omega $ is orthogonal to ${\mathcal{H}}_{n-p,n-q}(\Omega
 _{2}),$ and this is a necessary condition (see~\cite{Kohn73}).
 Because $h\in {\mathcal{H}}_{n-p,n-q}(\Omega _{2})$ implies
 that $h\in \Lambda _{n-p,n-q}(\bar \Omega _{2}),$ the scalar
 product ${\left\langle{\omega ,h}\right\rangle}$ is well defined
 and we have:\ \par 
\quad \quad \quad $\forall h\in {\mathcal{H}}_{n-p,n-q}(\Omega _{2}),\ {\left\langle{\omega
 ,h}\right\rangle}_{\Omega _{2}}={\left\langle{\bar \partial
 S,h}\right\rangle}_{\Omega _{2}}={\left\langle{S,\bar \partial
 ^{*}h}\right\rangle}_{\Omega _{2}}=0.$\ \par 
Hence $\omega $ is in the range of $\Delta _{\bar \partial },$
 so noting $N$ as usual the inverse of $\Delta _{\bar \partial
 },$ we get that $N$ is well defined on $\omega $ and Theorem
 4 in~\cite{BealsGreiStan87} gives that there is a $(p,q-1)$-current
 $g_{0}\in W_{p,q-1}^{1/2,r}(\Omega _{2})$ such that $\bar \partial
 g_{0}=\omega .$ Moreover, on any compact set $K\Subset \Omega
 _{2},$ we have $g_{0}\in W_{p,q-1}^{1,r}(K)$ because on $K$
 any vectors field can be extended to $\Omega _{2}$ as an admissible
 vectors field.\ \par 
\quad In particular we can choose $K:=\bar \Omega _{1},$ so we have
 that $g_{0}\in W_{p,q-1}^{1,r}(\bar \Omega _{1}).$\ \par 
\quad Now we have $\bar \partial (S-g_{0})=\omega -\bar \partial g_{0}=0$
 in $\Omega _{2}.$\ \par 
If $q=1,$ then $S-g_{0}$ is holomorphic in $\Omega _{2},$ hence
 ${\mathcal{C}}^{\infty }$ in $\bar \Omega _{1},$ so we have
 directly that $S\in W_{p,q-1}^{1,r}(\bar \Omega _{1}).$\ \par 
\quad Suppose now that $q\geq 2.$ Because $\Omega _{1}\backslash U$
 is not in general pseudo-convex even if $\Omega _{1}$ is, we
 cannot end the proof as in~\cite{Claurent15}.\ \par 
\quad So again we want to apply  Theorem 4 from ~\cite{BealsGreiStan87}
 to $\omega ':=S-g_{0}$ in $\Omega _{1}.$ We have to verify that
 $\omega '$ is orthogonal to ${\mathcal{H}}_{n-p,n+1-q}(\Omega _{1}).$\ \par 
But recall that in $\Omega _{2},\ g_{0}:=\bar \partial ^{*}N\omega
 $ hence, because $\forall h\in {\mathcal{H}}_{n-p,n+1-q}(\Omega
 _{1})\Rightarrow h\in \mathrm{d}\mathrm{o}\mathrm{m}(\bar \partial
 )\cap \mathrm{d}\mathrm{o}\mathrm{m}(\bar \partial ^{*})$ and
 $\bar \partial h=0$:\ \par 
\quad \quad \quad ${\left\langle{g_{0},h}\right\rangle}_{\Omega _{1}}={\left\langle{\bar
 \partial ^{*}N\omega ,h}\right\rangle}_{\Omega _{1}}={\left\langle{N\omega
 ,\bar \partial h}\right\rangle}_{\Omega _{1}}=0.$\ \par 
By assumption $S\perp {\mathcal{H}}_{n-p,n+1-q}(\Omega _{1})$
 hence $\omega '\perp {\mathcal{H}}_{n-p,n-q}(\Omega _{1}).$\ \par 
\ \par 
So again there is a $(p,q-2)$-current $g_{1}\in W_{p,q-2}^{1/2,r}(\Omega
 _{1})$ such that $\bar \partial g_{1}=S-g_{0}.$ And again, on
 any compact set $K\Subset \Omega _{1},$ we have $g_{1}\in W_{p,q-2}^{1,r}(K).$\
 \par 
\quad Let $\chi \in {\mathcal{C}}^{\infty }(X)$ such that $\chi =0$
 near the support  $W$ of $S,$ and $\chi =1$ in a neighborhood
 of $X\backslash U.$ Then the form $u:=g_{0}+\bar \partial (\chi
 g_{1})$ verifies:\ \par 
\quad \quad \quad $\bar \partial u=\bar \partial g_{0}=\omega $ because $\bar \partial
 ^{2}=0.$\ \par 
Now\ \par 
\quad \quad \quad $\bar \partial (\chi g_{1})=\bar \partial \chi \wedge g_{1}+\chi
 \bar \partial g_{1}=\bar \partial \chi \wedge g_{1}+\chi \bar
 \partial g_{1}=\bar \partial \chi \wedge g_{1}+\chi (S-g_{0}).$\ \par 
So we get, because $\bar \partial \chi =0$ and $S=0$ outside $U$\ \par 
\quad \quad \quad $\bar \partial (\chi g_{1})=-g_{0}$ outside $U$\ \par 
hence $u=0$ outside $U.$ Hence $u$ has its support in $U.$\ \par 
Now in $U$ we have $\chi =0$ near $W$ so\ \par 
\quad \quad \quad $\bar \partial (\chi g_{1})=\bar \partial \chi \wedge g_{1}-\chi g_{0}$\ \par 
and, because $\chi \in {\mathcal{C}}^{\infty }(X)$ and $g_{0},g_{1}\in
 W_{p,q-1}^{1,r}(K)$ for any compact $K$ in $\Omega _{1},$ we get\ \par 
\quad \quad \quad $u=\bar \partial \chi \wedge g_{1}-\chi g_{0}(1-\chi )g_{0}\in
 W_{p,q-1}^{1,r}(K)$ for any compact $K$ in $\Omega _{1}.$ In
 particular, because $u$ has its support in $U,$ we get that
 $u\in W_{p,q-1}^{1,r}(X).$\ \par 
\quad The origin of this method of control of the support is in section
 3.5, p. 9 of~\cite{AnGrauLr12}.\ \par 
\quad The proof is complete. $\hfill\blacksquare $\ \par 
\ \par 
Let us see the following example.\ \par 

\begin{exa}
~\label{SI0} There is a bounded open set with smooth boundary
 $\Omega $ in ${\mathbb{C}}^{n}$ and a compact set $K\subset
 \Omega $ such that there is no pseudo-convex set $D$ contained
 in $\Omega $ and containing $K.$
\end{exa}
\quad Proof.\ \par 
Take a bounded open cooking pot as $\Omega $ in ${\mathbb{C}}^{2}$
 and a compact one $K$ in $\Omega $ ( one can smoothed the boundaries)
 see Figure 1.\ \par 
\ \par 
\begin{figure}[h]
\begin{center}

\resizebox{14cm}{!}{\includegraphics{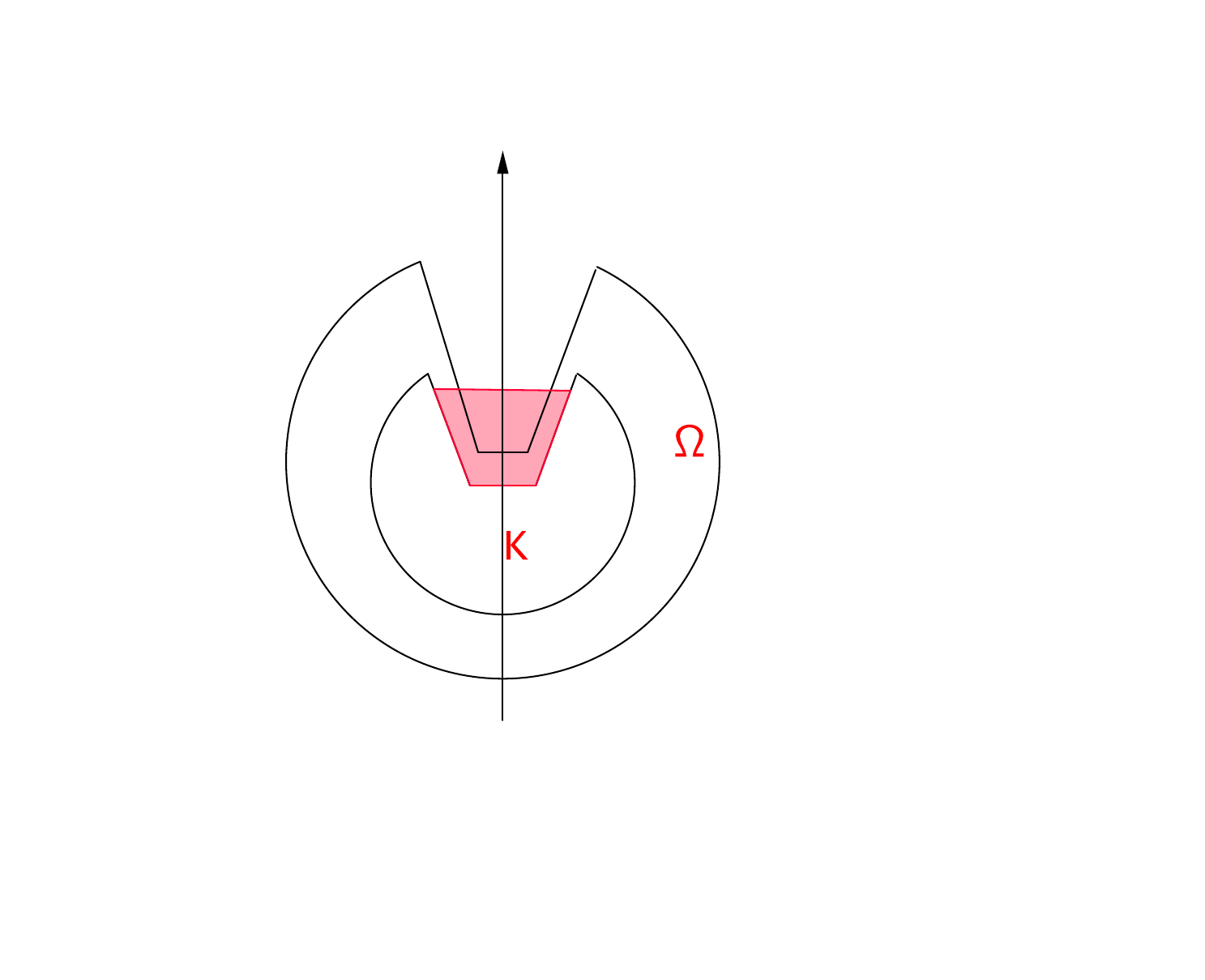}}
\caption{The cooking pot}
\end{center}

\end{figure}\ \par 
\quad Rotate the picture around the vertical axis in ${\mathbb{R}}^{4}={\mathbb{C}}^{2}$
 to get $\Omega .$\ \par 
Suppose that there is a pseudo-convex set $D$ in $\Omega $ and
 containing $K.$ Take any holomorphic function $h$ in $D.$ Then
 $h$ is holomorphic in a neighborhood of the boundary of  $K.$
 By the Cartan-Thullen Theorem $h$ extends in the red part, hence
 outside $D,$ so $D$ is not a domain of holomorphy hence is not
 pseudo-convex. $\hfill\blacksquare $\ \par 

\begin{rem}
Because of Example~\ref{SI0}, and the fact that an open set in
 ${\mathbb{C}}^{n}$ is a complex manifold, it seems difficult
 to get rid of the assumption that the support of $S$ must be
 in a pseudo-convex domain of $X.$\par 
Also the condition that the forms we want to solve this way be
 in the range of $\Delta _{\bar \partial }$ is necessary.
\end{rem}

\begin{rem}
Because the Theorem of Beals and all. is valid for domains having
 $Z(q)$ boundary, it is enough to suppose that $\Omega _{2}$
 and $\Omega _{1}$ be of type $Z(q)$ and $Z(q-1).$ The condition
 to belong to the range of $\Delta _{K}$ being the same as for
 the pseudo-convex case, by Theorem 3.2.2 p. 57 and the results
 at the beginning of p. 51 in~\cite{FolKoh72}.
\end{rem}
\ \par 
\quad The next corollary says that we can suppress the assumption $S\perp
 {\mathcal{H}}_{n-p,n+1-q}(\Omega _{1})$ provided that, on $\bar
 \Omega _{1},$ is defined a smooth strictly pluri-subharmonic function.\ \par 

\begin{cor}
Let $X$ be a complex manifold. Let $S$ be a $(p,q-1)$-current
 in $L_{p,q-1}^{r}(X).$  with compact support $W$ in $X.$ Suppose
 that $W\subset \Omega _{1}\Subset \Omega _{2},$ where $\Omega
 _{j},\ j=1,2$ are pseudo-convex open sets with smooth ${\mathcal{C}}^{\infty
 }$ boundary in $X$ and such that  there is a strictly pluri-subharmonic
 function  $\rho _{1}$ in ${\mathcal{C}}^{3}(\bar \Omega _{1}).$\par 
Moreover suppose that $\omega :=\bar \partial S$ is also in $L_{p,q}^{r}(X).$
 Let $U$ be any neighborhood of $W$ contained in $\Omega _{1}.$
 Then there is a $(p,q-1)$-current $u$ with compact support in
 $U$ such that $\bar \partial u=\omega $ and $u\in W_{p,q-1}^{1,r}(U).$
\end{cor}
\quad Proof.\ \par 
In fact we shall prove that, for $q\geq 1,\ {\mathcal{H}}_{p,q}(\Omega
 _{1})=\lbrace 0\rbrace .$ So let $h\in {\mathcal{H}}_{p,q}(\Omega
 _{1}).$ then $h\in \Lambda _{p,q}(\bar \Omega _{1}),\ \bar \partial
 h=\bar \partial ^{*}h=0.$ \ \par 
Because $h\in \Lambda _{p,q}(\bar \Omega _{1})$ and $\bar \Omega
 _{1}$ is compact, we have $h\in L_{p,q}^{2}(\Omega _{1}).$\ \par 
\quad If $\Omega _{1}$ is strongly pseudo-convex i.e. $\Omega _{1}:=\lbrace
 z\in X,\ \rho _{1}(z)<0\rbrace $ with $\rho _{1}\in {\mathcal{C}}^{2}(\bar
 \Omega )$ and the smallest eigenvalue of the form $\partial
 \bar \partial \rho _{1}$ is bounded below by $c_{\rho }>0$ by
 the continuity of $\partial \bar \partial \rho _{1}$ in $\bar
 \Omega ,$ we can apply Corollary~\ref{k52} with $r=2.$\ \par 
\quad With just the existence of $\rho _{1}$ as in the statement of
 the theorem, we can apply a well known $L^{2}$ Theorem of~\cite{Hormander73}:\
 \par 
\quad \quad \quad $\exists k\in L_{p,q-1}^{2}(\Omega _{1})::\bar \partial k=h.$\ \par 
So we have, because $h\in \mathrm{d}\mathrm{o}\mathrm{m}(\bar
 \partial ^{*})$ and $\bar \partial ^{*}h=0,$\ \par 
\quad \quad \quad ${\left\Vert{h}\right\Vert}^{2}={\left\langle{h,h}\right\rangle}={\left\langle{h,\bar
 \partial k}\right\rangle}={\left\langle{\bar \partial ^{*}h,k}\right\rangle}=0.$\
 \par 
The proof is complete. $\hfill\blacksquare $\ \par 

\begin{rem}
Because the Theorem 3.4.10 p. 145 in~\cite{Hormander65} is valid
 for domains having $Z(q)$ boundary (called $a_{q}$ in~\cite{Hormander65}),
 we have the same kind of corollary for these domains, provided
 that the defining function $\varphi _{j}$ for $\Omega _{j}$
 is defined in a neighborhood of $\bar \Omega _{j},$ verifies
 the condition $Z(q)$ outside of $\Omega _{j,c}:=\lbrace z\in
 \Omega _{j},\ \varphi (z)<c\rbrace $ for some $c>c_{0}$ and
 is exhausting in $\Omega _{j},\ j=1,2.$
\end{rem}
\quad So adding the results of Corollary~\ref{k53} and of Theorem~\ref{k45},
 we get:\ \par 

\begin{thm}
~\label{k55}Let $\Omega $ be a strictly pseudo-convex domain
 in a complex manifold and  $\omega $ be a $(p,q)$ form in $L^{r}(\Omega
 ),\ r>1$ with compact support in $\Omega .$ Suppose that $\omega
 $ is such that:\par 
\quad $\bullet $ if $1\leq q<n,\ \bar \partial \omega =0$; \par 
\quad $\bullet $ if $q=n,\ \forall V\subset \Omega ,\ \Supp \omega
 \subset V,\ \omega \perp {\mathcal{H}}_{n-p}(V)$\!\!\!\! .\par 
Then there is a $(p,q-1)$ form $u$ in $W^{1,r}(\Omega )$ with
 compact support in $\Omega $ such that $\bar \partial u=\omega
 $ as distributions and 	$\displaystyle {\left\Vert{u}\right\Vert}_{W^{1,r}(\Omega
 )}\leq C{\left\Vert{\omega }\right\Vert}_{L^{r}(\Omega )}.$
\end{thm}
\quad And the natural corollary:\ \par 

\begin{cor}
Let $X$ be a Stein manifold and  $\omega $ be a $(p,q)$ form
 in $L^{r}(X),\ r>1$ with compact support in $X.$ Suppose that
 $\omega $ is such that:\par 
\quad $\bullet $ if $1\leq q<n,\ \bar \partial \omega =0$; \par 
\quad $\bullet $ if $q=n,\ \forall V\subset X,\ \Supp \omega \subset
 V,\ \omega \perp {\mathcal{H}}_{n-p}(V)$\!\!\!\! .\par 
Then there is a $(p,q-1)$ form $u$ in $W^{1,r}(X)$ with compact
 support in $X$ such that $\bar \partial u=\omega $ as distributions
 and 	$\displaystyle {\left\Vert{u}\right\Vert}_{W^{1,r}(\Omega
 )}\leq C{\left\Vert{\omega }\right\Vert}_{L^{r}(\Omega )}.$
\end{cor}
\quad Proof.\ \par 
Because $X$ is a Stein manifold, for any compact set $K$ in $X$
 there is a relatively compact strictly pseudo-convex set $\Omega
 $ containing $K.$ So we can apply Theorem~\ref{k55}. $\hfill\blacksquare
 $\ \par 

\section{Hodge laplacian on riemannian manifolds.~\label{kF40}}
\quad A riemannian manifold $(M,g)$ is a real, smooth manifold $M$
 equipped with an inner product $g_{x}$ on the tangent space
 $T_{x}M$ at each point $x$ that varies smoothly from point to
 point in the sense that if $X$ and $Y$ are differentiable vector
 fields on $M,$ then $x\rightarrow g_{x}(X(x),Y(x))$ is a smooth
 function. The family $g_{x}$ of inner products is called a riemannian
 metric.\ \par 
\ \par 
\quad Let $X$ be a complete oriented riemannian manifold and $\Omega
 $ a relatively compact domain in $X.$\ \par 
We shall denote by $\Lambda ^{p}(\Omega )$ the set of ${\mathcal{C}}^{\infty
 }$ smooth $p$-forms in $\Omega $ and by $L^{r}_{p}(\Omega )$
 its closure in the Lebesgue space $L^{r}(\Omega )$ with respect
 to the riemannian volume measure $dm$ on $X.$\ \par 
\quad We shall take the following notation from the book by C. Voisin~\cite{Voisin02}.\
 \par 
\quad To a $p$-form $\alpha $ on $\Omega $ we associate its Hodge $*\
 (n-p)$-form $*\alpha .$ This gives us a pointwise scalar product
 and a pointwise modulus:\ \par 
\quad \quad \quad \begin{equation}  (\alpha ,\beta )dm:=\alpha \wedge *\bar \beta
 ;\ \ \left\vert{\alpha }\right\vert ^{2}dm:=\alpha \wedge {\overline{*\alpha
 ,}}\label{AG23}\end{equation}\ \par 
because $\displaystyle \alpha \wedge *\bar \beta $ is a $n$-form
 hence is a function time the volume form $dm.$\ \par 
\quad With the volume measure, we have a scalar product ${\left\langle{\alpha
 ,\beta }\right\rangle}$ on $p$-forms such that$\displaystyle
 \int_{\Omega }{\left\vert{\alpha }\right\vert ^{2}dm}<\infty
 .$ The link between these notions is given by~\cite[Lemme 5.8,
 p. 119]{Voisin02}:\ \par 
\quad \quad \quad \begin{equation} {\left\langle{ \alpha ,\beta }\right\rangle}=\int_{\Omega
 }{\alpha \wedge {\overline{*\beta }}}.\label{AG18}\end{equation}\ \par 
We shall define now $L_{p}^{r}(\Omega )$ to be the set of $p-$forms
 $\alpha $ defined on $\Omega $ such that\ \par 
\quad \quad \quad \quad 	$\ {\left\Vert{\alpha }\right\Vert}_{L_{p}^{r}(\Omega )}^{r}:=\int_{\Omega
 }{\left\vert{\alpha (z)}\right\vert ^{r}dm(z)}<\infty ,$\ \par 
where $\left\vert{\alpha }\right\vert $ is defined by~(\ref{AG23}).\ \par 
\quad As usual let $\displaystyle {\mathcal{D}}_{p}(\Omega )$ be the
 set of ${\mathcal{C}}^{\infty }\ p$-forms with compact support
 in $\Omega .$\ \par 
\quad On the manifold $M$ we have the exterior derivative $d$ on $p$-forms.
 To it we associate its \emph{formal adjoint} $d^{*}$ defined by:\ \par 
\quad \quad \quad $\forall u\in L^{r}_{p}(\Omega ),\ \forall \varphi \in {\mathcal{D}}_{p-1}(\Omega
 ),\ {\left\langle{d^{*}u,\varphi }\right\rangle}:={\left\langle{u,d\varphi
 }\right\rangle}.$\ \par 
Now we define the Hodge laplacian to be\ \par 
\quad \quad \quad $\Delta :=dd^{*}+d^{*}d.$\ \par 
This operator sends $p$-form to $p$-form and is essentially self-adjoint.
 In case $p=0,$ i.e. on functions, $\Delta $ is the usual Laplace-Beltrami
 operator on $M.$\ \par 
\quad We proved the following theorem~\cite[Theorem 1.1]{ellipticEq18},
 written here in the special case of the Hodge laplacian:\ \par 

\begin{thm}
~\label{kB2}Let $(M,g)$ be a ${\mathcal{C}}^{\infty }$ smooth
 {\bf compact} riemannian manifold without boundary. Let $\Delta
 :\Lambda _{p}\rightarrow \Lambda _{p}$ be the Hodge laplacian
 acting on the $p$-forms over $M.$ Let $\omega \in L^{r}_{p}(M)\cap
 \mathrm{(}\mathrm{k}\mathrm{e}\mathrm{r}\Delta \mathrm{)}^{\perp
 }$ with $r\in (1,\infty ).$ Then there is a bounded linear operator
 $S:L^{r}_{p}(M)\cap \mathrm{(}\mathrm{k}\mathrm{e}\mathrm{r}\Delta
 \mathrm{)}^{\perp }\rightarrow W^{2,r}_{p}(M)$ such that $\displaystyle
 \Delta S(\omega )=\omega $ on $M.$ So, with $u:=S\omega $ we
 get $\displaystyle \Delta u=\omega $ and $\displaystyle u\in
 W^{2,r}_{p}(M).$ Moreover we have $\displaystyle {\left\Vert{u}\right\Vert}_{W_{p}^{2,r}(M)}\leq
 c{\left\Vert{\omega }\right\Vert}_{L_{p}^{r}(M)}.$
\end{thm}
\ \par 
\quad We also proved the following theorem~\cite[Theorem 4.3, p. 14]{ellipticEq18},
 as a consequence of the Local Increasing Regularity Method.
 We just need to know here that the WMP is a weaker property
 than the Unique Continuation Property.\ \par 

\begin{thm}
Let $N$ be a smooth {\bf compact} riemannian manifold with smooth
 boundary $\partial N.$ Let $D:G\rightarrow G$ be an elliptic
 linear differential operator of order $m$ with ${\mathcal{C}}^{1}$
 coefficients acting on sections of a vector bundle $G:=(H,\pi
 ,M)$ on $N.$ Let $\displaystyle \omega \in L^{r}_{G}(N)$ be
 such a section. There is a $G$-section $u\in W^{m,r}_{G}(N),$
 such that $Du=\omega $ and ${\left\Vert{u}\right\Vert}_{W^{m,r}_{G}(N)}\leq
 c{\left\Vert{\omega }\right\Vert}_{L^{r}_{G}(N)},$ provided
 that the operator $D$ has the WMP for the $D$-harmonic $G$-forms.
\end{thm}
\quad Taking here $G:=\Lambda _{p}(N)$ the set of $p$-forms on $N,$
 and $D:=\Delta $ the Hodge laplacian, we have that  $\Delta
 $ verifies the Unique Continuation Property by a difficult result
  by N. Aronszajn, A. Krzywicki and J. Szarski~\cite{Aronszajn62}
 hence it has the WMP too.\ \par 
\quad So we get in this special case:\ \par 

\begin{thm}
~\label{k36}Let $N$ be a smooth {\bf compact} riemannian manifold
 with smooth boundary $\partial N.$ Let $\Delta $ be the Hodge
 Laplacian acting on $p$-forms on $N.$ Let $\displaystyle \omega
 \in L^{r}_{p}(N).$ There is a $p$-form $u\in W^{2,r}_{p}(N),$
 such that $\Delta u=\omega $ and ${\left\Vert{u}\right\Vert}_{W^{2,r}_{p}(N)}\leq
 c{\left\Vert{\omega }\right\Vert}_{L^{r}_{p}(N)}.$
\end{thm}
\quad This Theorem has the easy corollary:\ \par 

\begin{cor}
~\label{k37}Let $X$ be a complete smooth riemannian manifold
 without boundary. Let $\Omega $ a relatively compact domain
 in $X.$ Let $\omega \in L^{r}_{p}(\Omega ).$ There is a $p$-form
 $u\in W^{2,r}_{p}(\Omega ),$ such that $\Delta u=\omega $ and
 ${\left\Vert{u}\right\Vert}_{W^{2,r}_{p}(\Omega )}\leq c{\left\Vert{\omega
 }\right\Vert}_{L^{r}_{p}(\Omega )}.$
\end{cor}
\quad Proof.\ \par 
Put $\bar \Omega $ in a compact sub manifold $N$ of $X$ with
 a smooth boundary. Extend $\omega $ by $0$ outside $\Omega ,$
   then this extension $\displaystyle \tilde \omega $ is still
 in $\displaystyle L^{r}_{p}(N).$ We can apply Theorem~\ref{k36}
 to get a $p$-form $\tilde u\in W^{2,r}_{p}(N)$ such that $\Delta
 \tilde u=\tilde \omega .$ Now we let $u$ to be the restriction
 of $\tilde u$ to $\Omega .$ This ends the proof of the corollary.
 $\hfill\blacksquare $\ \par 

\begin{rem}
In the case where $\Omega $ is a bounded domain in ${\mathbb{R}}^{n},$
 to get this solution we just have to use the Newton kernel on
 $\omega $ and apply ~\cite[Theorem 9.9, p. 230]{GuilbargTrudinger98}.
 In the riemannian case we have to add a difficult result  by
 N. Aronszajn, A. Krzywicki and J. Szarski~\cite{Aronszajn62} to get the UCP.
\end{rem}

\section{Solution of the Poisson equation with compact support.~\label{kF42}}
\quad \quad 	Firts we shall study a duality between currents inspired by
 the Serre duality~\cite{Serre55}.\ \par 
\quad Because using Theorem~\ref{kB2}, the following results are easy,
 we shall assume from now on that $X$ is \emph{non compact.}\ \par 
\quad So let $X$ be an oriented non compact riemannian manifold of
 dimension $n.$ It has a volume form $dm$ and we denote also
 by $dm$ the associated volume measure on $X.$ We shall denote
 by $r'$ the conjugate exponent of $r\in (1,\infty ),\ \frac{1}{r}+\frac{1}{r'}=1.$\
 \par 

\subsection{Weighted $L^{r}$ spaces.}
\ \par 
\quad Let $\Omega $ be a domain in $X.$ \ \par 

\begin{lem}
~\label{AG22}Let $\eta >0$ be a weight. If $u$ is a $p$-current
 defined on $(n-p)$-forms $\alpha $ in $L^{r'}(\Omega ,\eta )$
 and such that\par 
\[ \forall \alpha \in L^{r'}_{n-p}(\Omega ,\eta ),\ \left\vert{{\left\langle{u,\
 *\alpha }\right\rangle}}\right\vert \leq C{\left\Vert{\alpha
 }\right\Vert}_{L^{r'}(\Omega ,\eta )},\]\par 
then $\ {\left\Vert{u}\right\Vert}_{L_{p}^{r}(\Omega ,\eta ^{1-r})}\leq C.$
\end{lem}
\quad \quad 	Proof.\ \par 
Set $\tilde \alpha :=\eta ^{1/r'}\alpha ;\ \tilde u:=\frac{1}{\eta
 ^{1/r'}}u$ then we have\ \par 
\quad \quad \quad $\displaystyle {\left\langle{u,*\alpha }\right\rangle}=\int_{\Omega
 }{u\wedge {\overline{\alpha }}}=\int_{\Omega }{\tilde u\wedge
 {\overline{\tilde \alpha }}}={\left\langle{\tilde u,*\tilde
 \alpha }\right\rangle}$\ \par 
and ${\left\Vert{\tilde \alpha }\right\Vert}_{L^{r'}(\Omega )}={\left\Vert{\alpha
 }\right\Vert}_{L^{r'}(\Omega ,\eta )}.$ \ \par 
\quad We notice that ${\left\Vert{\tilde \alpha }\right\Vert}_{L^{r'}(\Omega
 )}={\left\Vert{*\tilde \alpha }\right\Vert}_{L^{r'}(\Omega )}$
 because we have $(*\tilde \alpha ,*\tilde \alpha )dm=*\tilde
 \alpha \wedge {\overline{**\tilde \alpha }}$ but $**\tilde \alpha
 =(-1)^{p(n-p)}\tilde \alpha ,$ by ~\cite[Lemma 5.5]{Voisin02},
 hence, because $\displaystyle (*\tilde \alpha ,*\tilde \alpha
 )$ is positive,  $(*\tilde \alpha ,*\tilde \alpha )=\left\vert{\tilde
 \alpha }\right\vert ^{2}.$\ \par 
By use of the duality $L^{r}_{p}(\Omega )-L^{r'}_{n-p}(\Omega
 ),$ done in Lemma~\ref{AG25}, we get\ \par 
\quad \quad \quad $\displaystyle {\left\Vert{\tilde u}\right\Vert}_{L_{p}^{r}(\Omega
 )}=\sup _{\alpha \in L_{n-p}^{r'}(\Omega ),\ \alpha \neq 0}\frac{\left\vert{{\left\langle{\tilde
 u,*\tilde \alpha }\right\rangle}}\right\vert }{{\left\Vert{\tilde
 \alpha }\right\Vert}_{L^{r'}(\Omega )}}.$\ \par 
But\ \par 
\quad \quad \quad $\displaystyle {\left\Vert{\tilde u}\right\Vert}_{L_{p}^{r}(\Omega
 )}^{r}:=\int_{\Omega }{\left\vert{u}\right\vert ^{r}\eta ^{-\frac{r}{r'}}dm}=\int_{\Omega
 }{\left\vert{u}\right\vert ^{r}\eta ^{1-r}dm}={\left\Vert{u}\right\Vert}_{L^{r}(\Omega
 ,\eta ^{1-r})}^{r}.$\ \par 
So we get\ \par 
\quad \quad \quad {\it  }$\displaystyle {\left\Vert{u}\right\Vert}_{L_{p}^{r}(\Omega
 ,\eta ^{1-r})}=\sup _{*\alpha \in L_{p}^{r'}(\Omega ,\eta ),\
 \alpha \neq 0}\frac{\left\vert{{\left\langle{u,*\alpha }\right\rangle}}\right\vert
 }{{\left\Vert{\alpha }\right\Vert}_{L^{r'}(\Omega ,\eta )}}.$\ \par 
\quad The proof is complete. $\hfill\blacksquare $\ \par 
\ \par 
\quad Let ${\mathcal{H}}_{p}(\Omega )$ be the set of all $p$ harmonic
 forms, i.e. $h\in {\mathcal{H}}_{p}(\Omega )\iff \Delta h=0$
 in $\Omega .$\ \par 
In order to simplify notation, we note the pairing for $\alpha
 $ a $p$-form and $\beta $ a $(n-p)$-form by:\ \par 
\quad \quad \quad $\displaystyle \ll \alpha ,\beta \gg :=\int_{\Omega }{\alpha
 \wedge \beta }.$\ \par 
With this notation we also have ${\left\langle{\alpha ,\beta
 }\right\rangle}=\ll \alpha ,{\overline{*\beta }}\gg .$\ \par 
\ \par 

\begin{lem}
~\label{k34} We have $\Delta ({\overline{*u}})={\overline{*\Delta
 u}}.$ And $\displaystyle \ll \Delta \alpha ,\beta \gg =\ll \alpha
 ,\Delta \beta \gg $ provided that $\alpha $ or $\beta $ has
 compact support. Moreover we have\par 
\quad \quad \quad $\omega \in L^{r}_{p}(\Omega ),\ \omega \perp {\mathcal{H}}_{p}^{r'}(\Omega
 )\iff \omega \perp {\mathcal{H}}_{n-p}^{r'}(\Omega ).$\par 
with the suitable notion of orthogonality:\par 
\quad \quad \quad $\displaystyle \omega \in L^{r}_{p}(\Omega ),\ \omega \perp {\mathcal{H}}_{p}^{r'}(\Omega
 )\iff \forall h\in {\mathcal{H}}_{p}^{r'}(\Omega ),\ {\left\langle{\omega
 ,h}\right\rangle}=0$\par 
and\par 
\quad \quad \quad $\displaystyle \omega \in L^{r}_{p}(\Omega ),\ \omega \perp {\mathcal{H}}_{n-p}^{r'}(\Omega
 )\iff \forall h\in {\mathcal{H}}_{n-p}^{r'}(\Omega ),\ \ll \omega ,h\gg =0.$
\end{lem}
\quad Proof.\ \par 
We have $\Delta \varphi =dd^{*}\varphi +d^{*}d\varphi .$ The
 definition of $d^{*}$ in~\cite[Section 5.1.2, p. 118]{Voisin02} gives:\ \par 
\quad \quad \quad $d^{*}=(-1)^{p}*^{-1}d*$ on $\Lambda ^{p}.$\ \par 
We also have by~\cite[Lemme 5.5, p. 117]{Voisin02}:\ \par 
\quad \quad \quad $*^{2}=(-1)^{p(n-p)}$ on $\Lambda ^{p}.$\ \par 
These facts give:\ \par 
\quad \quad \quad $\displaystyle d(*\varphi )=**^{-1}d(*\varphi )=(-1)^{p}*d^{*}\varphi .$\ \par 
And, replacing the first $d^{*},$\ \par 
\quad \quad \quad $d^{*}d(*\varphi )=(-1)^{p}d^{*}*d^{*}\varphi =(-1)^{p}(-1)^{p}*^{-1}d**d^{*}\varphi
 =$\ \par 
\quad \quad \quad \quad \quad \quad \quad \quad \quad \quad \quad \quad \quad \quad \quad $\displaystyle =(-1)^{2p}(-1)^{2p(n-p)}*dd^{*}\varphi =*dd^{*}\varphi ,$\ \par 
because $\displaystyle *^{2}=(-1)^{p(n-p)}\Rightarrow *^{-1}=(-1)^{p(n-p)}*.$
 Hence $\displaystyle d^{*}d(*\varphi )=*dd^{*}\varphi .$\ \par 
\quad The same way we get $dd^{*}(*\varphi )=*d^{*}d\varphi .$ Because
 the laplacian is real the bar gets out.\ \par 
Now suppose that $\alpha $ has compact support we have:\ \par 
\quad \quad \quad $\ll \Delta \alpha ,\beta \gg ={\left\langle{\Delta \alpha ,{\overline{*\beta
 }}}\right\rangle}={\left\langle{\alpha ,\Delta ({\overline{*\beta
 }})}\right\rangle}={\left\langle{\alpha ,{\overline{*\Delta
 \beta }}}\right\rangle}=\ll \alpha ,\Delta \beta \gg ,$\ \par 
the second equality because $\Delta $ is essentially self-adjoint
 and the third one by the first part of this lemma.\ \par 
\quad For the "moreover", we have $h\in {\mathcal{H}}_{p}^{r'}(\Omega
 )\iff {\overline{*h}}\in {\mathcal{H}}_{n-p}^{r'}(\Omega )$
 because the first part of the lemma gives:\ \par 
\quad \quad \quad $\displaystyle \Delta ({\overline{*h}})={\overline{*\Delta h}}=0.$\ \par 
Now take $\displaystyle \omega \in L^{r}_{p}(\Omega )$ and $h\in
 {\mathcal{H}}_{p}^{r'}(\Omega )$ such that ${\left\langle{\omega
 ,h}\right\rangle}=0$ then\ \par 
\quad \quad \quad $\displaystyle 0={\left\langle{\omega ,h}\right\rangle}=\ll \omega
 ,{\overline{*h}}\gg $\ \par 
and the same for the converse, starting with  $\displaystyle
 h\in {\mathcal{H}}_{n-p}^{r'}(\Omega )$ and $\displaystyle \ll
 \omega ,h\gg =0$ we get $\displaystyle {\left\langle{\omega
 ,{\overline{*h}}}\right\rangle}=0.$\ \par 
The proof is complete. $\hfill\blacksquare $\ \par 
\ \par 
\quad Suppose that $\Omega $ is relatively compact in $X.$ Let $\omega
 \in L_{p}^{r}(\Omega )$ with compact support in $\Omega ,\ \omega
 \in L_{p}^{r,c}(\Omega ).$ \ \par 
\ \par 
Set the weight $\eta =\eta _{\epsilon }:={\11}_{\Omega _{1}}(z)+\epsilon
 {\11}_{\Omega \backslash \Omega _{1}}(z)$ for a fixed $\epsilon
 >0,$ with $\Supp \omega \subset \Omega _{1}\Subset \Omega .$\ \par 
\ \par 
\quad Let $\alpha \in L^{r'}_{p}(\Omega ,\eta ),$ with $r'$ conjugate
 to $r.$ Because $\epsilon >0$ we have $\alpha \in L^{r'}(\Omega
 ,\eta )\Rightarrow \alpha \in L^{r'}(\Omega ).$\ \par 
By Corollary~\ref{k37} we get:\ \par 
\quad \quad \quad \begin{equation}  \forall \alpha \in L^{r'}_{n-p}(\Omega ,\eta
 ),\ \exists \varphi \in W^{2,r'}_{n-p}(\Omega ),\ \Delta \varphi
 =\alpha ::{\left\Vert{\varphi }\right\Vert}_{W^{2,r'}(\Omega
 )}\lesssim {\left\Vert{\alpha }\right\Vert}_{L^{r'}(\Omega )}.\label{HD5}\end{equation}\
 \par 

\begin{lem}
Let $\omega \in L_{p}^{r}(\Omega _{1})\cap {\mathcal{H}}_{n-p}(\Omega
 _{1})^{\perp }$ with compact support in $\Omega _{1}$ and define:\par 
\quad \quad \quad $\displaystyle \forall \alpha \in L^{r'}_{n-p}(\Omega ,\eta ),\
 {\mathcal{L}}(\alpha ):=\ll \varphi ,\omega \gg ,$\par 
where $\varphi $ is a solution in $\Omega $ of~(\ref{HD5}).\par 
Then ${\mathcal{L}}$ is well defined and linear on $\displaystyle
 L^{r'}_{n-p}(\Omega ,\eta ).$
\end{lem}
\quad Proof.\ \par 
In order for ${\mathcal{L}}(\alpha )$ to be well defined, we
 need that if $\varphi '$ is another solution of $\Delta \varphi
 '=\alpha ,$ then $\ll \varphi -\varphi ',\omega \gg =0;$ hence
 we need that $\omega $ must be "orthogonal" to $(n-p)$-forms
 $\varphi $ such that $\Delta \varphi =0$ in $\Omega ,$ which
 is contained in our assumption.\ \par 
\quad Hence we have that ${\mathcal{L}}(\alpha )$ is well defined.\ \par 
The linearity of ${\mathcal{L}}$ is clear because if $\alpha
 =\alpha _{1}+\alpha _{2}$ take $\varphi _{j}::\Delta \varphi
 _{j}=\alpha _{j}$ then $\varphi :=\varphi _{1}+\varphi _{2}$
 implies $\displaystyle \Delta \varphi =\alpha _{1}+\alpha _{2}$ and\ \par 
\quad \quad \quad ${\mathcal{L}}(\alpha ):=\ll \varphi ,\omega \gg =\ll \varphi
 _{1},\omega \gg +\ll \varphi _{2},\omega \gg ={\mathcal{L}}(\alpha
 _{1})+{\mathcal{L}}(\alpha _{2}).$\ \par 
The same for $\lambda \alpha .$ The proof is complete. $\hfill\blacksquare
 $\ \par 
\ \par 
\quad By the H\"older inequalities done in Lemma~\ref{AG24} we get,
 because $\omega $ has its support in $\Omega _{1},$\ \par 
\quad \quad \quad $\left\vert{\ll \varphi ,\omega \gg }\right\vert =\left\vert{{\left\langle{\varphi
 ,{\overline{*\omega }}}\right\rangle}}\right\vert \leq {\left\Vert{\omega
 }\right\Vert}_{L^{r}(\Omega _{1})}{\left\Vert{\varphi }\right\Vert}_{L^{r'}(\Omega
 _{1})}.$\ \par 
Let $\alpha ,\varphi $ be as in~(\ref{HD5}), then\ \par 
\quad \quad \quad ${\left\Vert{\varphi }\right\Vert}_{L^{r'}(\Omega )}\leq {\left\Vert{\varphi
 }\right\Vert}_{W^{2,r'}(\Omega )}\leq C{\left\Vert{\alpha }\right\Vert}_{L^{r'}(\Omega
 )}.$\ \par 
But $\displaystyle {\left\Vert{\alpha }\right\Vert}_{L^{r'}(\Omega
 )}$ can be very big compared to $\displaystyle {\left\Vert{\alpha
 }\right\Vert}_{L^{r'}(\Omega _{1})}.$ So let $\psi $ be such
 that $\Delta \psi =\alpha $ in $\Omega _{1}$ and with ${\left\Vert{\psi
 }\right\Vert}_{W^{2,r'}(\Omega _{1})}\leq C{\left\Vert{\alpha
 }\right\Vert}_{L^{r'}(\Omega _{1})}.$ This is possible by Corollary~\ref{k37},
 $\bar \Omega _{1}$ being compact.\ \par 
\quad Then, because $\Delta (\varphi -\psi )=0$ in $\Omega _{1}$ and
 $\omega \perp {\mathcal{H}}_{n-p}(\Omega _{1}),$ we get\ \par 
\quad \quad \quad $\displaystyle {\mathcal{L}}(\alpha ):=\ll \varphi ,\omega \gg
 =\ll \psi ,\omega \gg .$\ \par 
Hence\ \par 
\quad \quad \quad $\ \left\vert{{\mathcal{L}}(\alpha )}\right\vert \leq {\left\Vert{\omega
 }\right\Vert}_{L^{r}(\Omega _{1})}{\left\Vert{\psi }\right\Vert}_{L^{r'}(\Omega
 _{1})}\leq C{\left\Vert{\omega }\right\Vert}_{L^{r}(\Omega _{1})}{\left\Vert{\alpha
 }\right\Vert}_{L^{r'}(\Omega _{1})}\leq C{\left\Vert{\omega
 }\right\Vert}_{L^{r}(\Omega _{1})}{\left\Vert{\alpha }\right\Vert}_{L^{r'}(\Omega
 ,\eta )},$        \ \par 
because $\eta _{\epsilon }=1$ on $\Omega _{1}\supset \Supp \omega
 ,$ hence ${\left\Vert{\alpha }\right\Vert}_{L^{r'}(\Omega _{1})}\leq
 {\left\Vert{\alpha }\right\Vert}_{L^{r'}(\Omega ,\eta )}.$\ \par 
\quad \quad 	So we have that the norm of ${\mathcal{L}}$ is bounded on $\displaystyle
 L_{n-p}^{r'}(\Omega ,\eta ).$ The bound of ${\mathcal{L}}$ 
 is $C{\left\Vert{\omega }\right\Vert}_{L^{r}(\Omega )}$ which
 is independent of $\eta $ hence of $\epsilon .$\ \par 
\quad \quad 	This means, by the definition of currents, that there is a $p$
 current $u$ which represents the form ${\mathcal{L}}$: $\displaystyle
 {\mathcal{L}}(\alpha )=\ll \alpha ,u\gg .$ So if $\alpha :=\Delta
 \varphi $ with $\varphi \in {\mathcal{C}}^{\infty }$ with compact
 support in $\Omega ,$ we get\ \par 
\quad \quad \quad $\displaystyle \ll \varphi ,\omega \gg ={\mathcal{L}}(\alpha
 )=\ll \alpha ,u\gg =\ll \Delta \varphi ,u\gg .$\ \par 
Now we use Lemma~\ref{k34} to get $\displaystyle \ll \varphi
 ,\omega \gg =\ll \varphi ,\Delta u\gg $ and we have $\Delta
 u=\omega $ in the sense of distributions. \ \par 
\quad Moreover we have\ \par 
\quad \quad \quad $\displaystyle \sup _{\alpha \in L^{r'}(\Omega ,\eta ),\ {\left\Vert{\alpha
 }\right\Vert}=1}\ \left\vert{\ll \alpha ,u\gg }\right\vert \leq
 C{\left\Vert{\omega }\right\Vert}_{L^{r}(\Omega )}$\ \par 
and by Lemma~\ref{AG22} with the weight $\eta ,$ this implies\ \par 
\quad \quad \quad \quad 	$\ {\left\Vert{u}\right\Vert}_{L^{r}(\Omega ,\eta ^{1-r})}\leq
 C{\left\Vert{\omega }\right\Vert}_{L^{r}(\Omega )}.$\ \par 
\quad So we proved\ \par 

\begin{prop}
~\label{k35}Let $\Omega _{1}\Subset \Omega $ and $\omega \in
 L^{r}(\Omega _{1})$ with compact support in $\Omega _{1}$ and
 such that $\omega \perp {\mathcal{H}}_{n-p}(\Omega _{1}).$ Let
 also $\eta =\eta _{\epsilon }:={\11}_{\Omega _{1}}(z)+\epsilon
 {\11}_{\Omega \backslash \Omega _{1}}(z).$ Then there is a $p$
 form $u\in L^{r}(\Omega ,\eta ^{1-r})$ such that $\displaystyle
 \Delta u=\omega $ and $\displaystyle {\left\Vert{u}\right\Vert}_{L^{r}(\Omega
 ,\eta ^{1-r})}\leq C{\left\Vert{\omega }\right\Vert}_{L^{r}(\Omega )}.$
\end{prop}
\quad Now we are in position to prove:\ \par 

\begin{thm}
~\label{k26}Let $X$ be a complete oriented riemannian manifold.
 Let $\Omega $ be a relatively compact domain in $X$ and $\Omega
 _{1}\Subset \Omega .$ Let $\omega \in L^{r}_{p}(\Omega _{1})$
 with compact support in $\Omega _{1}$ and such that $\omega
 \perp {\mathcal{H}}_{n-p}(\Omega _{1}).$ Then there is a $p$-form
 $u\in L^{r}_{p}(\Omega )$ with compact support in $\Omega _{1}$
 such that $\Delta u=\omega $ as distributions and 	$\displaystyle
 {\left\Vert{u}\right\Vert}_{L^{r}_{p}(\Omega )}\leq C{\left\Vert{\omega
 }\right\Vert}_{L^{r}_{p}(\Omega _{1})}.$
\end{thm}
\quad \quad 	Proof.\ \par 
	For $\epsilon >0$ with $\eta _{\epsilon }(z):={\11}_{\Omega
 _{1}}(z)+\epsilon {\11}_{\Omega \backslash \Omega _{1}}(z),$
 let $u_{\epsilon }\in L^{r}(\Omega ,\eta _{\epsilon }^{1-r})$
 be the solution given by Proposition~\ref{k35}, then\ \par 
\quad \quad \quad \quad 	$\displaystyle \ {\left\Vert{u_{\epsilon }}\right\Vert}_{L^{r}(\Omega
 ,\eta _{\epsilon }^{1-r})}^{r}\leq \int_{\Omega }{\left\vert{u_{\epsilon
 }}\right\vert ^{r}\eta ^{1-r}dm}\leq C^{r}{\left\Vert{\omega
 }\right\Vert}_{L^{r}(\Omega )}^{r}.$\ \par 
Replacing $\eta $ by its value we get\ \par 
\quad \quad \quad \quad 	$\displaystyle \ \int_{\Omega _{1}}{\left\vert{u_{\epsilon }}\right\vert
 ^{r}dm}+\ \int_{\Omega \backslash \Omega _{1}}{\left\vert{u_{\epsilon
 }}\right\vert ^{r}\epsilon ^{1-r}dm}\leq C^{r}{\left\Vert{\omega
 }\right\Vert}_{L^{r}(\Omega )}^{r}$ $\displaystyle \Rightarrow
 \int_{\Omega \backslash \Omega _{1}}{\left\vert{u_{\epsilon
 }}\right\vert ^{r}\epsilon ^{1-r}dm}\leq C^{r}{\left\Vert{\omega
 }\right\Vert}_{L^{r}(\Omega )}^{r}$\ \par 
hence\ \par 
\quad \quad \quad \quad \quad  $\displaystyle \ \int_{\Omega \backslash \Omega _{1}}{\left\vert{u_{\epsilon
 }}\right\vert ^{r}dm}\leq C^{r}\epsilon ^{r-1}{\left\Vert{\omega
 }\right\Vert}_{L^{r}(\Omega )}^{r}.$\ \par 
Because $C$ and the norm of $\omega $ are independent of $\epsilon
 ,$ we have that $\ {\left\Vert{u_{\epsilon }}\right\Vert}_{L^{r}(\Omega
 )}$ is uniformly bounded and $r>1$ implies that $\displaystyle
 L_{p}^{r}(\Omega )$ is a dual by Lemma~\ref{AG25}, hence there
 is a sub-sequence $\lbrace u_{\epsilon _{k}}\rbrace _{k\in {\mathbb{N}}}$
 of $\lbrace u_{\epsilon }\rbrace $ which converges weakly to
 a $p$-form $u$ in $L_{p}^{r}(\Omega ),$ when $\epsilon _{k}\rightarrow
 0,$ still with $\ {\left\Vert{u}\right\Vert}_{L_{p}^{r}(\Omega
 )}\leq C{\left\Vert{\omega }\right\Vert}_{L_{p}^{r}(\Omega )}.$
 Let us note $\displaystyle u_{k}:=u_{\epsilon _{k}}.$\ \par 
\quad To see that this form $u$ is $0\ a.e.$ on $\Omega \backslash
 \Omega _{1}$ let us write the weak convergence:\ \par 
\quad \quad \quad $\displaystyle \forall \alpha \in L_{p}^{r'}(\Omega ),\ {\left\langle{u_{k},\alpha
 }\right\rangle}=\int_{\Omega }{u_{k}\wedge {\overline{*\alpha
 }}}\underset{k\rightarrow \infty }{\rightarrow }{\left\langle{u,\alpha
 }\right\rangle}=\int_{\Omega }{u\wedge {\overline{*\alpha }}}.$\ \par 
As usual take $\displaystyle \alpha :=\frac{u}{\left\vert{u}\right\vert
 }{\11}_{E}$ where $E:=\lbrace \left\vert{u}\right\vert >0\rbrace
 \cap (\Omega \backslash \Omega _{1})$ then we get\ \par 
\quad \quad \quad $\displaystyle \int_{\Omega }{u\wedge {\overline{*\alpha }}}=\int_{E}{\left\vert{u}\right\vert
 dm}=\lim _{k\rightarrow \infty }\int_{\Omega }{u_{k}\wedge {\overline{*\alpha
 }}}=\lim _{k\rightarrow \infty }\int_{E}{\frac{u_{k}\wedge {\overline{*u}}}{\left\vert{u}\right\vert
 }}.$\ \par 
Now we have by H\"older inequalities:\ \par 
\quad \quad \quad $\displaystyle \left\vert{\int_{E}{\frac{u_{k}\wedge {\overline{*u}}}{\left\vert{u}\right\vert
 }}}\right\vert \leq {\left\Vert{u_{k}}\right\Vert}_{L^{r}(E)}{\left\Vert{{\11}_{E}}\right\Vert}_{L^{r'}(E)}.$\
 \par 
But\ \par 
\quad \quad \quad \quad 	$\displaystyle \ {\left\Vert{u_{k}}\right\Vert}_{L^{r}(E)}^{r}\leq
 \int_{\Omega \backslash \Omega _{1}}{\left\vert{u_{k}}\right\vert
 ^{r}dm}\leq (\epsilon _{k})^{r-1}C{\left\Vert{\omega }\right\Vert}_{L^{r}(\Omega
 )}\underset{k\rightarrow \infty }{\rightarrow }0$\ \par 
and ${\left\Vert{{\11}_{E}}\right\Vert}_{L^{r'}(E)}=(m(E))^{1/r'}.$\ \par 
Hence\ \par 
\quad \quad \quad \quad 	$\displaystyle \ \left\vert{\ \int_{E}{\left\vert{u}\right\vert
 dm}}\right\vert =\lim _{k\rightarrow \infty }\int_{E}{\frac{u_{k}\wedge
 {\overline{*u}}}{\left\vert{u}\right\vert }}\leq \ \lim _{k\rightarrow
 \infty }C^{r}(m(E))^{1/r'}(\epsilon _{k})^{r-1}{\left\Vert{\omega
 }\right\Vert}_{L^{r}(\Omega )}^{r}=0,$\ \par 
so $\int_{E}{\left\vert{u}\right\vert dm}=0$ which implies $m(E)=0$
 because on $E,\ \left\vert{u}\right\vert >0.$\ \par 
\quad Hence we get that the form $u$ is $0\ a.e.$ on $\Omega \backslash
 \Omega _{1}.$\ \par 
\quad \quad 	So we proved\ \par 
\quad \quad \quad \quad 	\begin{equation}  \forall \varphi \in {\mathcal{D}}_{n-p}(\Omega
 ),\ \ll \varphi ,\omega \gg =\ll \Delta \varphi ,u_{k}\gg \underset{k\rightarrow
 \infty }{\rightarrow }\ll \Delta \varphi ,u\gg \Rightarrow \ll
 \Delta \varphi ,u\gg =\ll \varphi ,\omega \gg \label{k27}\end{equation}\ \par 
hence again by use of Lemma~\ref{k34} we get $\Delta u=\omega
 $ in the sense of distributions. $\hfill\blacksquare $\ \par 

\begin{lem}
Let $X$ be a complete riemannian manifold. Let $\Omega $ be a
 relatively compact domain in $X$ and $\Omega _{1}\Subset \Omega
 .$ Let $u\in L^{r}_{p}(\Omega )$ such that $\Delta u\in L^{r}_{p}(\Omega
 ).$ Then we have the interior elliptic regularity:\par 
\quad \quad \quad $\displaystyle {\left\Vert{u}\right\Vert}_{W^{2,r}_{p}(\Omega
 _{1})}\leq C({\left\Vert{\Delta u}\right\Vert}_{L^{r}_{p}(\Omega
 )}+{\left\Vert{u}\right\Vert}_{L^{r}_{p}(\Omega )}).$
\end{lem}
\quad Proof.\ \par 
The interior elliptic inequalities~\cite{ellipticEq18}, Theorem
 3.4, valid in the complete riemannian manifold $M$ give that
 for any $x\in M,$ there is a ball $B_{x}:=B(x,R)$ and a smaller
 ball $B'_{x}$ relatively compact in $B_{x},$ such that:\ \par 
\quad \quad \quad \begin{equation} {\left\Vert{ u}\right\Vert}_{W^{2,r}(B'x)}\leq
 c_{1}{\left\Vert{\Delta u}\right\Vert}_{L^{r}(B_{x})}+c_{2}R(x)^{-2}{\left\Vert{u}\right\Vert}_{L^{r}(B_{x})}.\label{k61}\end{equation}\
 \par 
Moreover the constants $c_{j},\ j=1,2,$ are independent of the
 radius $R(x)$ of the ball $B_{x}.$\ \par 
Because $\bar \Omega _{1}$ is compact in $\Omega ,$ there is
 a $\delta >0$ such that:\ \par 
\quad \quad \quad $\bigcup_{x\in \Omega _{1}}{B(x,\delta )}\subset \Omega .$\ \par 
For all $x\in \bar \Omega _{1},$ choose $R'(x)=\min (\delta ,R(x))$
 for the $R(x)$ given in~(\ref{k61}).\ \par 
\quad We cover the compact set $\bar \Omega _{1}$ by a finite set of
 balls $B'_{x_{j}}$ associated to $B(x_{j},R'(x)).$ So we get,
 by~(\ref{k61}),\ \par 
\quad \quad \quad ${\left\Vert{u}\right\Vert}_{W^{2,r}(\Omega _{1})}\leq \sum_{j=1}^{N}{}{\left\Vert{u}\right\Vert}_{W^{2,r}(B'_{x_{j}})}\leq
 c_{1}\sum_{j=1}^{N}{}{\left\Vert{\Delta u}\right\Vert}_{L^{r}(B_{x_{j}})}+c_{2}\sum_{j=1}^{N}{}R'(x_{j})^{-2}{\left\Vert{u}\right\Vert}_{L^{r}(B_{x_{j}})}.$\
 \par 
Set $c:=\max _{j=1,..,N}R'(x_{j})^{-2},$ we get:\ \par 
\quad \quad \quad ${\left\Vert{u}\right\Vert}_{W^{2,r}(\Omega _{1})}\leq c_{1}\sum_{j=1}^{N}{}{\left\Vert{\Delta
 u}\right\Vert}_{L^{r}(B_{x_{j}})}+cc_{2}\sum_{j=1}^{N}{}{\left\Vert{u}\right\Vert}_{L^{r}(B_{x_{j}})}.$\
 \par 
Now we have\ \par 
\quad \quad \quad $\sum_{j=1}^{N}{}{\left\Vert{f}\right\Vert}_{L^{r}(B_{x_{j}})}=\sum_{j=1}^{N}{}\int_{B_{x_{j}}}{\left\vert{f}\right\vert
 ^{r}}dm\leq \sum_{j=1}^{N}{}\int_{\Omega }{\left\vert{f}\right\vert
 ^{r}}dm\leq N\int_{\Omega }{\left\vert{f}\right\vert ^{r}}dm,$\ \par 
because, by the choice of $R'(x)$ we have that $B(x_{j},R')\subset
 \Omega .$\ \par 
\quad Applying this with $f=u$ and $f=\Delta u$ we get\ \par 
\quad \quad \quad ${\left\Vert{u}\right\Vert}_{W^{2,r}(\Omega _{1})}\leq c_{1}N{\left\Vert{\Delta
 u}\right\Vert}_{L^{r}(\Omega )}+cc_{2}N{\left\Vert{u}\right\Vert}_{L^{r}(\Omega
 )}.$\ \par 
The proof is complete. $\hfill\blacksquare $\ \par 
\ \par 
\quad This lemma allows the better estimates:\ \par 

\begin{cor}
~\label{k29}Let $X$ be a complete oriented riemannian manifold.
 Let $\Omega $ be a relatively compact domain in $X$ and $\Omega
 _{1}\Subset \Omega .$ Let $\omega \in L^{r}_{p}(\Omega _{1})$
 with compact support in $\Omega _{1}$ and such that $\omega
 \perp {\mathcal{H}}_{n-p}(\Omega _{1}).$ Then there is a $p$-form
 $u\in W^{2,r}_{p}(\Omega )$ with compact support in $\Omega
 _{1}$ such that $\Delta u=\omega $ as distributions and 	$\displaystyle
 {\left\Vert{u}\right\Vert}_{W^{2,r}_{p}(\Omega _{1})}\leq C{\left\Vert{\omega
 }\right\Vert}_{L^{r}_{p}(\Omega _{1})}.$
\end{cor}
\quad Proof.\ \par 
We can apply Theorem~\ref{k26} so we have a $p$-form $u\in L^{r}_{p}(\Omega
 )$ with compact support in $\Omega _{1}$ such that $\Delta u=\omega
 $ as distributions and ${\left\Vert{u}\right\Vert}_{L^{r}_{p}(\Omega
 )}\lesssim {\left\Vert{\omega }\right\Vert}_{L^{r}(\Omega )}.$\ \par 
\quad Now we apply the interior elliptic regularity with $\Delta u=\omega $:\ \par 
\quad \quad \quad $\displaystyle {\left\Vert{u}\right\Vert}_{W^{2,r}_{p}(\Omega
 _{1})}\leq C({\left\Vert{\omega }\right\Vert}_{L^{r}_{p}(\Omega
 )}+{\left\Vert{u}\right\Vert}_{L^{r}_{p}(\Omega )}).$\ \par 
But ${\left\Vert{u}\right\Vert}_{L^{r}_{p}(\Omega )}\lesssim
 {\left\Vert{\omega }\right\Vert}_{L^{r}(\Omega )}$ so we get\ \par 
\quad \quad \quad $\displaystyle {\left\Vert{u}\right\Vert}_{W^{2,r}_{p}(\Omega
 _{1})}\leq C{\left\Vert{\omega }\right\Vert}_{L^{r}_{p}(\Omega )}.$\ \par 
Because $\omega $ has compact support in $\Omega _{1},$ we get
 $\displaystyle {\left\Vert{u}\right\Vert}_{W^{2,r}_{p}(\Omega
 _{1})}\leq C{\left\Vert{\omega }\right\Vert}_{L^{r}_{p}(\Omega _{1})}.$\ \par 
The proof is complete. $\hfill\blacksquare $\ \par 

\begin{rem}
~\label{k28}The condition of orthogonality to ${\mathcal{H}}_{p}(\Omega
 _{1})$ is necessary: suppose there is a $p$-current $u\in W^{2,r}_{p}(\Omega
 )$ such that $\Delta u=\omega $ and $u$ with compact support
 in $\Omega ,$ then if $h\in {\mathcal{H}}_{n-p}(\Omega ),$ we have\par 
\quad \quad  $h\in {\mathcal{H}}_{n-p}(\Omega ),\ \ll \omega ,h\gg =\ll \Delta
 u,h\gg =\ll u,\Delta h\gg =0,$\par 
because $u$ is compactly supported.
\end{rem}

\section{K\"ahler manifold and Kohn laplacian.~\label{kF41}}
\quad A K\"ahler manifold is a complex manifold $X$ with a Hermitian
 metric $h$  whose associated 2-form $\kappa $ is closed. In
 more detail, $h$ gives a positive definite Hermitian form on
 the tangent space $T_{x}$ at each point $x$ of $X,$ and the
 2-form $\kappa $ is defined by\ \par 
\quad \quad \quad $\kappa (u,v):=\Re h(iu,v)$\ \par 
for tangent vectors $u$ and $v$ (where $i$ is the complex number
 ${\sqrt{-1}}$ ). For a K\"ahler manifold $X,$ the K\"ahler form
 $\kappa $ is a real closed (1,1)-form. A K\"ahler manifold can
 also be viewed as a Riemannian manifold, with the Riemannian
 metric $g$ defined by\ \par 
\quad \quad \quad $\displaystyle g(u,v):=\Re h(u,v).$\ \par 
\ \par 
\quad On $X$ the $(p,q)$-forms are defined and so is the $\bar \partial
 $ operator. The Hodge $*$ operator is also defined, see C. Voisin~\cite[Section
 5.1.4, p. 121]{Voisin02}.\ \par 
\quad Recall the $\bar \partial $ (or Kohn) laplacian, acting from
 $(p,q)$-forms to $(p,q)$-forms is:\ \par 
\quad \quad \quad $\Delta _{\bar \partial }f:=(\bar \partial \bar \partial ^{*}+\bar
 \partial ^{*}\bar \partial )f,$\ \par 
where $\bar \partial ^{*}$ is the formal adjoint to $\bar \partial
 ,$ i.e.\ \par 
\quad \quad \quad $\forall \varphi \in {\mathcal{D}}_{p,q-1},\ \forall u\in L^{r}_{p,q},\
 {\left\langle{\bar \partial ^{*}u,\ \varphi }\right\rangle}:={\left\langle{u,\bar
 \partial \varphi }\right\rangle}.$\ \par 
\quad The space ${\mathcal{H}}^{r}_{q}(\Omega ):=\lbrace h\in L_{q}^{r}(\Omega
 )::\Delta h=0\rbrace $ is the space of harmonic $q$-forms in
 the set $\Omega .$\ \par 
\quad Because $X$ is a complex manifold, it is canonically oriented
 and we also note $dm$ the volume $(n,n)$ form on $X.$\ \par 
\ \par 
\quad Now our aim is to prove Theorem ~\ref{kF39}:\ \par 

\begin{thm}
Let $(X,\omega )$ be a complete K\"ahler manifold. Let $\Omega
 $ be a relatively compact  domain in $X.$ Let $\displaystyle
 \omega \in L_{p,q}^{r}(\Omega ),\ \bar \partial \omega =0$ in
 $\Omega $ and $\omega $ compactly supported in $\Omega .$ Suppose
 moreover that $\displaystyle \omega \perp {\mathcal{H}}^{r'}_{2n-p-q}(\Omega
 ).$\par 
\quad Then there is a $\displaystyle u\in W_{p,q-1}^{1,r}(\Omega )$
 with compact support in $\Omega $ and such that $\bar \partial u=\omega .$
\end{thm}
\quad Proof.\ \par 
Let us see $X$ as a riemannian manifold, then we can apply Corollary~\ref{k29}
 to get the existence of a $\displaystyle \tilde v\in W_{p+q}^{2,r}(\Omega
 )$ such that $\Delta \tilde v=\omega $ and $\tilde v$ compactly
 supported in $\Omega .$\ \par 
\quad By use of Theorem~\ref{k60} we get that $\Delta _{\bar \partial
 }\tilde v=\frac{1}{2}\omega .$ So, setting $v:=\frac{1}{2}\tilde
 v$ we get:\ \par 
\quad \quad \quad $\displaystyle v\in W_{p+q}^{2,r,c}(\Omega )::\Delta _{\bar \partial
 }v=\omega .$\ \par 
Now we have\ \par 
\quad \quad \quad \begin{equation}  \Delta _{\bar \partial }v=(\bar \partial \bar
 \partial ^{*}+\bar \partial ^{*}\bar \partial )v=\omega \label{k31}\end{equation}\
 \par 
this implies, taking $\bar \partial $ on both sides,\ \par 
\quad \quad \quad $\bar \partial \bar \partial ^{*}\bar \partial v=\bar \partial
 \omega =0,$\ \par 
because $\bar \partial ^{2}=0.$ Then\ \par 
\quad \quad \quad $0={\left\langle{\bar \partial \bar \partial ^{*}\bar \partial
 v,\bar \partial v}\right\rangle}={\left\langle{\bar \partial
 ^{*}\bar \partial v,\bar \partial ^{*}\bar \partial v}\right\rangle}={\left\Vert{\bar
 \partial ^{*}\bar \partial v}\right\Vert}_{L^{2}(\Omega )}^{2}$\ \par 
because $v$ being compactly supported in $\Omega ,$ so is $\bar
 \partial v,$ and we can shift the $\bar \partial $ operator
 on the right hand side.\ \par 
\quad From~(\ref{k31}) we get $\displaystyle \ \bar \partial \bar \partial
 ^{*}v=\omega ,$ because $\displaystyle \bar \partial ^{*}\bar
 \partial v=0.$  Now we set $u:=\bar \partial ^{*}v$ and we get
 a $u$ with support in $\Omega ,$ such that:\ \par 
\quad \quad \quad $u\in W_{p,q-1}^{1,r}(\Omega ),\ \bar \partial u=\omega ,$\ \par 
because $\bar \partial ^{*}$ is a first order differential operator
 and $\displaystyle v\in W_{p+q}^{2,r}(\Omega )$ with support
 in $\Omega .$\ \par 
\quad The proof is complete. $\hfill\blacksquare $\ \par 

\begin{rem}
1) In the case of bounded open sets in ${\mathbb{C}}^{n}$ and
 for the $L^{2}$ theory, this  idea to use the usual laplacian
 to get estimates for the $\bar \partial $ equation was already
 done in the nice book by E. Straube~\cite[Section 2.9]{Straube10}.\par 
\quad 2) This method improves the regularity of the solution: from
 $\displaystyle L^{r,c}_{p,q-1}(\Omega )$ to $\displaystyle W_{p,q-1}^{1,r,c}(\Omega
 ).$ The price is that $\displaystyle \omega \perp {\mathcal{H}}_{2n-p-q}(\Omega
 )$ but there is no pseudo-convexity condition on $\Omega .$
\end{rem}

\section{Appendix.}
\quad For the reader's convenience we shall prove certainly known results
 on the duality $L^{r}-L^{r'}$ for $(p,q)$-forms in a complex manifold.\ \par 
\ \par 
Recall we have a pointwise scalar product and a pointwise modulus:\ \par 
\quad \quad \quad $\displaystyle (\alpha ,\beta )dm:=\alpha \wedge {\overline{*\beta
 }};\ \ \left\vert{\alpha }\right\vert ^{2}dm:=\alpha \wedge
 {\overline{*\alpha .}}$\ \par 
By the Cauchy-Schwarz inequality for scalar product we get:\ \par 
\quad \quad \quad $\displaystyle \forall x\in X,\ \left\vert{(\alpha ,\beta )(x)}\right\vert
 \leq \left\vert{\alpha (x)}\right\vert \left\vert{\beta (x)}\right\vert
 .$\ \par 
\quad This gives H\"older inequalities for $(p,q)$-forms:\ \par 

\begin{lem}
~\label{AG24}(H\"older inequalities) Let $\alpha \in L^{r}_{p,q}(\Omega
 )$ and $\displaystyle \beta \in L^{r'}_{p,q}(\Omega ).$ We have\par 
\quad \quad \quad $\displaystyle \left\vert{{\left\langle{\alpha ,\beta }\right\rangle}}\right\vert
 \leq {\left\Vert{\alpha }\right\Vert}_{L^{r}(\Omega )}{\left\Vert{\beta
 }\right\Vert}_{L^{r'}(\Omega )}.$
\end{lem}
\quad Proof.\ \par 
We start with ${\left\langle{\alpha ,\beta }\right\rangle}=\int_{\Omega
 }{(\alpha ,\beta )(x)dm(x)}$ hence\ \par 
\quad \quad \quad $\displaystyle \left\vert{{\left\langle{\alpha ,\beta }\right\rangle}}\right\vert
 \leq \int_{\Omega }{\left\vert{(\alpha ,\beta )(x)}\right\vert
 dm}\leq \int_{\Omega }{\left\vert{\alpha (x)}\right\vert \left\vert{\beta
 (x)}\right\vert dm(x)}.$\ \par 
By the usual H\"older inequalities for functions we get\ \par 
\quad \quad \quad $\displaystyle \int_{\Omega }{\left\vert{\alpha (x)}\right\vert
 \left\vert{\beta (x)}\right\vert dm(x)}\leq {\left({\int_{\Omega
 }{\left\vert{\alpha (x)}\right\vert ^{r}dm}}\right)}^{1/r}{\left({\int_{\Omega
 }{\left\vert{\beta (x)}\right\vert ^{r'}dm}}\right)}^{1/r'}$\ \par 
which ends the proof of the lemma. $\hfill\blacksquare $\ \par 

\begin{lem}
~\label{AG21}Let $\alpha \in L_{p,q}^{r}(\Omega )$ then\par 
\quad \quad \quad $\displaystyle {\left\Vert{\alpha }\right\Vert}_{L_{p,q}^{r}(\Omega
 )}=\sup _{\beta \in L_{p,q}^{r'}(\Omega ),\ \beta \neq 0}\frac{\left\vert{{\left\langle{\alpha
 ,\beta }\right\rangle}}\right\vert }{{\left\Vert{\beta }\right\Vert}_{L^{r'}(\Omega
 )}}.$
\end{lem}
\quad Proof.\ \par 
We choose $\beta :=\alpha \left\vert{\alpha }\right\vert ^{r-2},$ then:\ \par 
\quad \quad \quad $\displaystyle \left\vert{\beta }\right\vert ^{r'}=\left\vert{\alpha
 }\right\vert ^{r'(r-1)}=\left\vert{\alpha }\right\vert ^{r}\Rightarrow
 {\left\Vert{\beta }\right\Vert}_{L^{r'}(\Omega )}^{r'}={\left\Vert{\alpha
 }\right\Vert}_{L^{r}(\Omega )}^{r}.$\ \par 
Hence\ \par 
\quad \quad \quad $\displaystyle {\left\langle{\alpha ,\beta }\right\rangle}={\left\langle{\alpha
 ,\alpha \left\vert{\alpha }\right\vert ^{r-2}}\right\rangle}=\int_{\Omega
 }{(\alpha ,\alpha )\left\vert{\alpha }\right\vert ^{r-2}dm}={\left\Vert{\alpha
 }\right\Vert}_{L^{r}(\Omega )}^{r}.$\ \par 
On the other hand we have\ \par 
\quad \quad \quad $\displaystyle {\left\Vert{\beta }\right\Vert}_{L^{r'}(\Omega
 )}={\left\Vert{\alpha }\right\Vert}_{L^{r}(\Omega )}^{r/r'}={\left\Vert{\alpha
 }\right\Vert}_{L^{r}(\Omega )}^{r-1},$\ \par 
so\ \par 
\quad \quad \quad $\displaystyle {\left\Vert{\alpha }\right\Vert}_{L^{r}(\Omega
 )}{\times}{\left\Vert{\beta }\right\Vert}_{L^{r'}(\Omega )}={\left\Vert{\alpha
 }\right\Vert}_{L^{r}(\Omega )}^{r}={\left\langle{\alpha ,\beta
 }\right\rangle}.$\ \par 
Hence\ \par 
\quad \quad \quad $\displaystyle {\left\Vert{\alpha }\right\Vert}_{L^{r}(\Omega
 )}=\frac{\left\vert{{\left\langle{\alpha ,\beta }\right\rangle}}\right\vert
 }{{\left\Vert{\beta }\right\Vert}_{L^{r'}(\Omega )}}.$\ \par 
A fortiori for any choice of $\beta $:\ \par 
\quad \quad \quad $\displaystyle {\left\Vert{\alpha }\right\Vert}_{L^{r}(\Omega
 )}\leq \sup _{\beta \in L^{r'}(\Omega )}\frac{\left\vert{{\left\langle{\alpha
 ,\beta }\right\rangle}}\right\vert }{{\left\Vert{\beta }\right\Vert}_{L^{r'}(\Omega
 )}}.$\ \par 
To prove the other direction, we use the H\"older inequalities,
 Lemma~\ref{AG24}:\ \par 
\quad \quad \quad $\displaystyle \forall \beta \in L_{p,q}^{r'}(\Omega ),\ \frac{\left\vert{{\left\langle{\alpha
 ,\beta }\right\rangle}}\right\vert }{{\left\Vert{\beta }\right\Vert}_{L^{r'}(\Omega
 )}}\leq {\left\Vert{\alpha }\right\Vert}_{L^{r}(\Omega )}.$\ \par 
The proof is complete. $\hfill\blacksquare $\ \par 
\ \par 
\quad Now we are in position to state:\ \par 

\begin{lem}
~\label{AG25}The dual space of the Banach space $\displaystyle
 L_{p,q}^{r}(\Omega )$ is $\displaystyle L_{n-p,n-q}^{r'}(\Omega ).$
\end{lem}
\quad Proof.\ \par 
Suppose first that $u\in L_{n-p,n-q}^{r'}(\Omega ).$ Then consider:\ \par 
\quad \quad \quad $\displaystyle \forall \alpha \in L_{p,q}^{r}(\Omega ),\ {\mathcal{L}}(\alpha
 ):=\int_{\Omega }{\alpha \wedge u}={\left\langle{\alpha ,{\overline{*u}}}\right\rangle}.$\
 \par 
This is a linear form on $\displaystyle L_{p,q}^{r}(\Omega )$
 and its norm, by definition, is\ \par 
\quad \quad \quad $\displaystyle {\left\Vert{{\mathcal{L}}}\right\Vert}=\sup _{\alpha
 \in L^{r}(\Omega )}\frac{\left\vert{{\left\langle{\alpha ,{\overline{*u}}}\right\rangle}}\right\vert
 }{{\left\Vert{\alpha }\right\Vert}_{L^{r}(\Omega )}}.$\ \par 
By use of Lemma~\ref{AG21} we get\ \par 
\quad \quad \quad $\displaystyle {\left\Vert{{\mathcal{L}}}\right\Vert}={\left\Vert{{\overline{*u}}}\right\Vert}_{L^{r'}_{p,q}(\Omega
 )}={\left\Vert{u}\right\Vert}_{L^{r'}_{n-p,n-q}(\Omega )}.$\ \par 
So we have ${\left({L_{p,q}^{r}(\Omega )}\right)}^{*}\supset
 L_{n-p,n-q}^{r'}(\Omega )$ with the same norm.\ \par 
\ \par 
\quad Conversely take a continuous linear form ${\mathcal{L}}$ on $\displaystyle
 L_{p,q}^{r}(\Omega ).$ We have, again by definition, that:\ \par 
\quad \quad \quad $\displaystyle {\left\Vert{{\mathcal{L}}}\right\Vert}=\sup _{\alpha
 \in L^{r}(\Omega )}\frac{\left\vert{{\mathcal{L}}(\alpha )}\right\vert
 }{{\left\Vert{\alpha }\right\Vert}_{L^{r}(\Omega )}}.$\ \par 
Because ${\mathcal{D}}_{p,q}(\Omega )\subset L_{p,q}^{r}(\Omega
 ),$  ${\mathcal{L}}$ is a continuous linear form on $\displaystyle
 {\mathcal{D}}_{p,q}(\Omega ),$ hence, by definition, ${\mathcal{L}}$
 can be represented by a $(n-p,n-q)$-current $u.$ So we have:\ \par 
\quad \quad \quad $\displaystyle \forall \alpha \in {\mathcal{D}}_{p,q}(\Omega
 ),\ {\mathcal{L}}(\alpha ):=\int_{\Omega }{\alpha \wedge u}={\left\langle{\alpha
 ,{\overline{*u}}}\right\rangle}.$\ \par 
Moreover we have, by Lemma~\ref{AG21},\ \par 
\quad \quad \quad $\displaystyle {\left\Vert{{\mathcal{L}}}\right\Vert}=\sup _{\alpha
 \in {\mathcal{D}}_{p,q}(\Omega )}\frac{\left\vert{{\left\langle{\alpha
 ,*\bar u}\right\rangle}}\right\vert }{{\left\Vert{\alpha }\right\Vert}_{L^{r}(\Omega
 )}}={\left\Vert{*u}\right\Vert}_{L^{r'}(\Omega )}$\ \par 
because $\displaystyle {\mathcal{D}}_{p,q}(\Omega )$ is dense
 in $\displaystyle L_{p,q}^{r}(\Omega ).$ So we proved\ \par 
\quad \quad \quad $\displaystyle {\left({L_{p,q}^{r}(\Omega )}\right)}^{*}\subset
 L_{n-p,n-q}^{r'}(\Omega )$ with the same norm.\ \par 
The proof is complete. $\hfill\blacksquare $\ \par 
\ \par 

\bibliographystyle{/usr/local/texlive/2017/texmf-dist/bibtex/bst/base/apalike}

\end{document}